\documentclass{amsart}

\usepackage{newtxtext, newtxmath}

\usepackage{latexsym, bm}
\usepackage{amsmath, amsfonts, amsthm}
 
\usepackage{amssymb}
\usepackage{cases}
\usepackage{enumerate, paralist}%asparaitem, inparaitem, compactitem, or --enum
\usepackage{appendix}

\usepackage[colorlinks=true]{hyperref}
\usepackage[alphabetic, msc-links, abbrev, nobysame]{amsrefs}
\usepackage{titletoc}%titletoc

\input{ThmAndAutoref}

\newcommand{\supp}{\mathrm{ supp }}

\DeclareMathOperator*{\id}{\mathrm{id}}

\allowdisplaybreaks

\numberwithin{equation}{section}

\begin{document}

\title[Machado--Bishop theorem]{The Machado--Bishop theorem in the uniform topology}

\author{Deliang Chen}
\address{College of Mathematics and Physics, Wenzhou University, Wenzhou 325035, People's Republic of China}
\email{chernde@wzu.edu.cn}

\thanks{The author is grateful to the anonymous referee for careful reading and helpful comments which greatly improved the presentation. The research is supported by the National Natural Science Foundation of China (No. 12101461).}

%    General info


\begin{abstract}
The Machado--Bishop theorem for weighted vector-valued functions vanishing at infinity has been extensively studied. In this paper, we give an analogue of Machado's distance formula for bounded weighted vector-valued functions. A number of applications are given; in particular, some types of the Bishop--Stone--Weierstrass theorem for bounded vector-valued continuous spaces in the uniform topology are discussed; the splitting of $C(I \times J, X \otimes Y)$ as the closure of $C(I, X) \otimes C(J, Y)$ in different senses and the extension of continuous vector-valued functions are studied.
\end{abstract}

\maketitle

%titletoc

\titlecontents{section}[0pt]{\vspace{0\baselineskip}\bfseries}
{\thecontentslabel\quad}{}%
{\hspace{0em}\titlerule*[10pt]{$\cdot$}\contentspage}

\titlecontents{subsection}[1em]{\vspace{0\baselineskip}}
{\thecontentslabel\quad}{}%
{\hspace{0em}\titlerule*[10pt]{$\cdot$}\contentspage}

\titlecontents{subsubsection}[2em]{\vspace{0\baselineskip}}
{\thecontentslabel\quad}{}%
{\hspace{0em}\titlerule*[10pt]{$\cdot$}\contentspage}

\setcounter{tocdepth}{2}

\section{Introduction}

Let $\Omega$ be a compact Hausdorff space. The celebrated Stone--Weierstrass theorem states that for a subalgebra $A$ of the complex-valued continuous functions space $C(\Omega, \mathbb{C})$, if $A$ is self-adjoint, contains a function vanishing nowhere and separates points on $\Omega$, then $A$ is dense in $C(\Omega, \mathbb{C})$ in the uniform topology. A number of generalizations were made. Among them, Bishop \cite{Bis61} realized that $f \in A$ where $A$ is a closed unital subalgebra of $C(\Omega, \mathbb{C})$, if $f|_E \in A|_E$ for some ``simple'' subsets $E \subset \Omega$ (i.e., the $A$-antisymmetric sets); the conditions of the Stone--Weierstrass theorem on $A$ mean that the $A$-antisymmetric sets are single points. In the real case (i.e., $A \subset C(\Omega, \mathbb{R})$), ``$f|_E \in A|_E$ for all $A$-antisymmetric sets $E$'' reads that ``$A$ separates $f$'' (i.e., if $f(x) \neq f(y)$ then there is $g\in A$ such that $g(x) \neq g(y)$). Later, Machado \cite{Mac77} gave a strong version of this theorem; that is, he established a distance formula,
\[
d_{\Omega}(f, A) = d_{E}(f, A)
\]
for some $A$-antisymmetric subset $E \subset \Omega$, where 
\[
d_{E}(f, A) = \inf_{g \in A}\max_{x \in E}\{ |f(x)-g(x)| \}.
\]

When $\Omega$ is \emph{non-compact}, such results are still possible if the uniform topology is replaced by a smaller one. For example, Meyer \cite{Mey67} used pointwise convergence (i.e. the point-open topology), and Arens \cite{Are49} used the topology of uniform convergence on compact subsets (i.e. the compact-open topology). Moreover, such results seem to hold if the space $C(\Omega, \mathbb{C})$ is taken by a smaller one; for instance, Buck \cite{Buc58} and Glicksberg \cite{Gli63} used the strict topology for bounded continuous functions space $C_b(\Omega, \mathbb{C})$, and de Branges \cite{dBra59} used the uniform topology for $C_0(\Omega, \mathbb{C})$ (the continuous functions vanishing at infinity). The Stone--Weierstrass type theorems also hold for \emph{vector-valued} functions, that is, the complex field $\mathbb{C}$ could be replaced by Banach spaces, or more generally, locally convex topological vector spaces $X$; Pe\l czy\'{n}ski \cite{Pczy57} might be the first one to investigate such generalization when $\Omega$ is compact.

In his seminal article \cite{Nac65}, Nachbin introduced the weighted space to deal with the uniform approximation over a non-compact space, and particularly generalized the Stone--Weierstrass theorem for the weighted spaces. By choosing different weighted families $V$ (i.e. Nachbin families), the topologies in weighted spaces $CV_0(\Omega, \mathbb{C})$ (weighted functions vanishing at infinity) include the topologies mentioned above (see \autoref{exa:c0}). 
Prolla \cite{Pro71} and Summers \cite{Sum71} gave a Bishop's version of the Stone--Weierstrass theorem in vector-valued weighted spaces $CV_0(\Omega, X)$ (see also \cite{PM73}). In \cite{Che18f}, Chen generalized the Machado theorem in weighted spaces $CV_0(\Omega, X)$. We refer the readers to see the monographs \cite{Pro77} and \cite{BP20} for more generalized results and important applications. 

In \cite{Hew47}, Hewitt gave a counter-example to show that the Stone--Weierstrass theorem fails to be true for the bounded real continuous function space $C_b(\Omega, \mathbb{R})$ in the \emph{uniform topology} if $\Omega$ is \emph{non-compact}; the right situation so that such theorem is valid in all completely regular spaces $\Omega$ is that ``$A$ separates points of $\Omega$'' should be replaced by that ``$A$ separates disjoint \emph{zero-sets} of $\Omega$'' (see also \cite{Nel68}). In the Bishop's version, when $A$ is a linear subspace of $C_b(\Omega, \mathbb{R})$, Blasco and Molt\'{o} \cite{BM83} considered $f \in \overline{A}$ where the closure is taken in the uniform topology, if ``$A$ \emph{S-separates} disjoint Lebesgue sets of $f$''. Garrido and Montalvo \cite{GM93} discussed the situation when $A$ is a subalgebra of $C_b(\Omega, \mathbb{R})$; in this case, that ``$A$ \emph{S-separates} disjoint Lebesgue sets of $f$'' could be weekended as that ``$A$ \emph{separates} disjoint Lebesgue sets of $f$''. By using the notion of multiplier, Bustamante and Montalvo \cite{BM00} obtained more general results for the uniform approximation for modules of $C_b(\Omega, \mathbb{R})$.

In this paper, we study the Machado--Bishop theorem for bounded weighted spaces $CV_b(\Omega, X)$ (see \autoref{machado} and \autoref{thm:ap}) and particularly $C_b(\Omega, X)$ and $C_b(\Omega, \mathbb{C})$ in the uniform topology. Unlike in the weighted space $CV_0(\Omega, X)$, one cannot expect that Machado's distance formula (see e.g. \cite[Theorem 2.8]{Che18f} or \autoref{cor:mac}) holds in $CV_b(\Omega, X)$ (and so $C_b(\Omega, \mathbb{R})$ in the uniform topology), due to Hewitt's counter-example \cite{Hew47}. Instead, we find a formula like 
\[
d_{\Omega}(f, A) = \inf_{E \in \mathcal{A}} d_{E}(f, A)
\]
where $\mathcal{A}$ is some ``simple'' $z$-filter (see e.g. \cite[2.2]{GJ60}) of $\Omega$; see \autoref{machado} in the general setting. There does not seem to be a formula in this form in the existing literature. We use this formula to discuss some types of Bishop and Stone--Weierstrass theorems for bounded (weighted) vector-valued continuous spaces in the uniform topology; see \autoref{sec:uap} and \autoref{sec:SW}.

\begin{enumerate}[(a)]
	\item In the existing literature, there are fewer results describing the \emph{uniform} closure of modules of $C_b(\Omega, \mathbb{C})$ when $\Omega$ is non-compact. If a subalgebra $A$ of $C_b(\Omega, \mathbb{C})$ is self-adjoint, then by the corresponding results for the real-valued functions case (e.g. \cite{Hew47, Nel68, GM93, BM00}), one might use that ``$A$ separates disjoint zero-sets of $\Omega$'' to describe the uniform closure of $A$, or that ``$A$ S-separates disjoint Lebesgue sets of $f$'' to determine whether $f$ belongs to the uniform closure of $A$. However, if $A$ is not self-adjoint, some new results or methods are needed. We give an example (see \autoref{exa:rudin}) to show how the distance formula we obtained can be applied to this situation.
	
	\item It seems that, except $CV_0(\Omega, X)$, types of Bishop and Stone--Weierstrass theorems for vector-valued functions space $C_b(\Omega, X)$ (and so $CV_b(\Omega, X)$) are rarely discussed in the existing literature. We find that, in general, the Hewitt's version of the Stone--Weierstrass theorem does not hold for $C_b(\Omega, X)$ (see \autoref{exa:nonSW}). However, if the locally convex topological vector space $X$ has the Heine--Borel property (i.e., every closed and bounded subset of $X$ is compact) (resp. $X$ is metrizable), then such result holds in the uniform topology (resp. the pseudocompact-open topology); see \autoref{thm:SW}.
	
	\item Consider the problem whether $f$ belongs to the uniform closure of $W$, a space of $C_b(\Omega, X)$, where $f \in C_b(\Omega, X)$. Such problem is related with ``localization'' (see e.g. \cite{Nac65, Pro77}) and is extensively studied in $CV_0(\Omega, X)$.
	Hern\'{a}ndez-Mu\~{n}oz (see \cite[Theorem 4]{Her94}) considered such uniform approximation for vector-valued functions in a special case. 
	It should be mentioned that Kleinst\"{u}ck \cite{Kle75} might be the first one to break new ground for such problem in $CV_b(\Omega, X)$; but the results are stated in terms of the Stone--\v{C}ech compactification and do not establish a relationship with ``separability'', which seemly leads to inconvenient applications in practical problems.
	Instead, we obtain a vector-valued version of the Bishop theorem due to Bustamante and Montalvo (see \cite[Theorem 5.4]{BM00}) when $f(\Omega)$ is relatively compact; see \autoref{thm:ap} in the general setting. Note that if $f(\Omega)$ is not relatively compact, then such result fails to be true (see \autoref{exa:compact}). 
	
	\item In \autoref{sec:tensor}, as an application of our main results, we study the splitting of $C(I \times J, X \otimes Y)$ as the closure of $C(I, X) \otimes C(J, Y)$ in different senses, which is related with the approximation property for function spaces (see e.g. \cite{Bie81}).
	
	\item In \autoref{sec:extension}, as another application, we study the extension of $f \in C(S, X)$ (with $f(S)$ relatively compact) to a function in $C(\Omega, X)$; in particular, \autoref{thm:ex} is a generalization of Mr\'{o}wka \cite[Theorem 4.11]{Mro68} (see also Blair \cite[Theorem 3.2]{Bla81}) to the vector-valued case.
\end{enumerate}

\section{Preliminaries and main result}

\subsection{Notations}

Throughout this paper, the Hausdorff property of a topological space is not assumed unless where
mentioned. We use the following notations.
\begin{enumerate}[$\bullet$]
	\item $\mathbb{R}_+ = \{x \in \mathbb{R}: x \geq 0\}$.
	\item $\mathbb{C}\cup\{\infty\}$: the one-point compactification of the complex field $\mathbb{C}$.
	\item $\Omega$: a topological space.
	\item $\overline{S}$: the closure of $S$ if $S$ is a subset of a topological space.
	\item $X$: a locally convex topological vector space with a family of seminorms denoted by $\mathfrak{A}$; if $\mathfrak{A}$ consists of one seminorm, the seminorm will be written as $|~|$.
	\item $F(\Omega, X)$: the space of all functions from $\Omega$ into $X$.
	\item $F_b(\Omega, X)$: the space of all bounded functions from $\Omega$ into $X$, that is, $f \in F_b(\Omega, X)$ if and only if for every $p \in \mathfrak{A}$, $\sup _{x \in \Omega}p(f(x)) < \infty$.
	\item $F_0(\Omega, X)$: the space of all functions that vanish at infinity, that is, $f \in F_0(\Omega, X)$ if and only if for each $p \in \mathfrak{A}$ and $\epsilon > 0$, the closure of $\{x \in \Omega: p(f(x)) \geq \epsilon\}$ is a compact subset of $\Omega$.
	\item $C(\Omega, X)$: the space of all continuous functions from $\Omega$ into $X$.
	\item $C_b(\Omega, X) = C(\Omega, X) \cap F_b(\Omega, X)$: the space of all bounded continuous functions from $\Omega$ into $X$.
	\item $ \supp f $: the support of $ f: \Omega \to X $ defined by $ \supp f = \overline{f^{-1}(X \setminus \{0\})} $.
	\item $ f|_S $: the restriction of $ f: \Omega \to X $ to $ S $ ($ \subset \Omega $). 
	\item $f(S) = \{f(x) : x \in S \}$: the range of $ f|_S $.
\end{enumerate}

For $p \in \mathfrak{A}$, $v \in F(\Omega, \mathbb{R}_+)$, $f \in F(\Omega, X)$, $F \subset \Omega$, define
\[
|f|_{v, p, F} = \sup\{v(x)p(f(x)): x \in F\}, ~|f|_{v, p} = |f|_{v, p, \Omega},
\]
and for $W \subset F(\Omega, X)$, 
\[
d_{v, p, F}(f, W) = \inf \{|f - g|_{v, p, F}: g \in W \};
\]
write
\[
f - W = \{ f - g : g \in W \}.
\]
For brevity, if $\mathfrak{A}$ consists of one seminorm or $v = 1$, we write $d_{F}(f, W)$ (or $d_{p, F}(f, W)$, $d_{v, F}(f, W)$) instead of $d_{v, p, F}(f, W)$.

\subsection{Weighted space}
If $V \subset F(\Omega, \mathbb{R}_+)$, define the (bounded) \emph{weighted space} (with respect to $V$) by 
\[
FV_{b}(\Omega, X) = \{ f \in F(\Omega, X) : vf \in F_{b}(\Omega, X) \text{ for all } v \in V \},
\]
and the (little) weighted space (with respect to $V$) by
\[
FV_{0}(\Omega, X) = \{ f \in F(\Omega, X) : vf \in F_{0}(\Omega, X) \text{ for all } v \in V \}.
\] 
For $W \subset F(\Omega, X)$, write 
\[
FV_{b}W = FV_{b}(\Omega, X) \cap W, ~FV_{0}W = FV_{0}(\Omega, X) \cap W,
\]
and 
\[
CV_{b}(\Omega, X) = FV_{b}C(\Omega, X), ~CV_{0}(\Omega, X) = FV_{0}C(\Omega, X).
\]

Let $V$ be a \emph{Nachbin family} on $\Omega$, i.e., a set of upper semicontinuous functions of $ \Omega \to \mathbb{R}_+ $ with the following two additional conditions on $V$:
\begin{enumerate}[(i)]
	\item for all $ v_1, v_2 \in V $, there are $ \lambda \geq 0 $ and $ w \in V $ such that $ v_1 \leq \lambda w $, $ v_2 \leq \lambda w $ (point-wise);
	\item for each $ x \in \Omega $, there is $ v \in V $ such that $ v(x) > 0 $.
\end{enumerate}
If $ W $ is a linear subspace of $ F(\Omega, X) $, then $ FV_b W $ is locally convex if condition (i) is satisfied and Hausdorff if condition (ii) holds, where a local subbase at $ 0 $ can be given by $ \{ f \in FV_b W : |f|_{v, p} < \epsilon \} $, $ v \in V $, $ p \in \mathfrak{A} $, and $ \epsilon > 0 $; the same assertion holds for $FV_0W$ and particularly $CV_{b}(\Omega, X)$, $CV_{0}(\Omega, X)$, etc. For more details about weighted spaces $CV_{b}(\Omega, X)$ and $CV_{0}(\Omega, X)$, see e.g. \cite{Nac65,Pro71, PM73, Bie73, Pro77, BP20}.

Some examples of $CV_{b}(\Omega, X)$ are the following. Let $ \chi_{K} $ be the characteristic function of $ K $, i.e., $ \chi_{K} (x) = 1 $ if $ x \in K $, and $ \chi_{K} (x) = 0 $ otherwise.

\begin{exa}\label{exa:bounded}
	\begin{enumerate}[(a)]
		\item Let $ V = \{ 1 \} $. Then $ CV_b (\Omega, X) = C_b(\Omega, X) $. In this case, the topology is the uniform topology. % and $CV_0 (\Omega, X) = C_0(\Omega, X)$
		\item \label{it:ps} Let $ V = \{ \chi_{K}: K \text{ is a closed pseudocompact subset of } \Omega \} $. Then $ CV_b (\Omega, X) = C(\Omega, X) $. The topology in this case is the \emph{pseudocompact-open} topology of $ C(\Omega, X) $. This topology is stronger than the compact-open topology. See \cite{KG06} for more details about the pseudocompact-open topology.
		\item Let $ V = \{ \chi_{K}: K \text{ is a closed pseudocompact and $C^*$-embedded subset of } \Omega \} $. Then $ CV_b (\Omega, X) = C(\Omega, X) $. Nel \cite{Nel68} called this topology of $ C(\Omega, X) $ the $\mathscr{L}$-topology. Clearly, this topology is stronger than the compact-open topology but weaker than the pseudocompact-open topology.
		\item \label{it:bd} A subset $S$ of $\Omega$ is called \emph{bounded} if for every $f \in C(\Omega, \mathbb{R})$, $f|_S$ is bounded (see also \cite{San18} for a survey); particularly, all the pseudocompact subsets of $\Omega$ are bounded. Let $ V = \{ \chi_{K}: K \text{ is a closed bounded subset of } \Omega \} $. Then $ CV_b (\Omega, X) = C(\Omega, X) $. This topology is stronger than the pseudocompact-open topology. Briefly, we call it the \emph{bounded-open} topology. See \cite[Section 2 Example 1 (iii)]{Bie73}.
	\end{enumerate}
\end{exa}

The topologies in (b) (c) (d) would be the same e.g. in all normal spaces (as all the closed subsets of normal spaces are $C$-embedded) and would be equal to the compact-open topology in metric spaces (as pseudocompact subsets of metric spaces are compact (see e.g. \cite[Chapter XI Section 3]{Dug66} or \cite[Proposition 1.1.13]{ACIT18})).  For some basic properties of pseudocompact, see e.g. \cite{ACIT18}.

For some examples of $CV_{0}(\Omega, X)$, see e.g. \citelist{\cite{Bie73}*{Section 2} \cite{Pro77}*{Section 5.1} \cite{BP20}*{Section 2.1}}, etc. In particular, by choosing different Nachbin families $V$, $CV_{0}(\Omega, X)$ includes different topologies of continuous functions spaces:

\begin{exa}\label{exa:c0}
	\begin{enumerate}[(a)]
		\item the point-open topology of $C(\Omega, X)$ where $V = \{ \chi_{K}: K \text{ is a finite subset of } \Omega \}$ when $X$ is a $T_1$-space;
		\item the compact-open topology of $C(\Omega, X)$ where $V = \{ \chi_{K}: K \text{ is a closed compact subset of } \Omega \}$; 
		\item the uniform topology of $C_0(\Omega, X)$ where $ V = \{ 1 \} $; 
		\item the strict topology of $C_{b}(\Omega, X)$ where $V = C_0(\Omega, \mathbb{R}_+)$ if $\Omega$ is locally compact Hausdorff \cite{Buc58} or $V = \text{ all the upper semicontinuous functions of } F_{0}(\Omega, \mathbb{R}_+)$ if $\Omega$ is completely regular Hausdorff \cite{Sum71}; 
		\item the inductive limit topology of $C_{c}(\Omega, X)$ (continuous functions with compact supports) where $V = C(\Omega, \mathbb{R}_+)$ and $X$ is locally compact Hausdorff and countable at infinity \cite[Page 192]{Bie73}.
	\end{enumerate}
\end{exa}

\subsection{$z$-filter}

A subset $S$ of $\Omega$ is called a \emph{zero-set} if there is $h \in C(\Omega, \mathbb{R})$ such that $S = h^{-1}(0)$. In particular, if $a \in \mathbb{R}$ and $h \in C(\Omega, \mathbb{R})$, then the \emph{Lebesgue sets} of $h$ defined by
\[
L_{a}(h) = \{x \in \Omega: h(x) \leq a\}, ~ L^{a}(h) = \{x \in \Omega: h(x) \geq a\},
\]
are zero-sets. The set of all zero-sets of $\Omega$ is denoted by $Z(\Omega)$. 
\begin{defi} [See {\cite[2.2]{GJ60}}]
	A nonempty collection $\mathcal{F} \subset Z(\Omega)$ is called a \emph{$z$-filter} on $\Omega$ if
	\begin{enumerate}[(i)]
		\item $\emptyset \notin \mathcal{F}$;
		\item if $Z_1, Z_2 \in \mathcal{F}$, then $Z_1 \cap Z_2 \in \mathcal{F}$; and
		\item if $Z \in \mathcal{F}$, $Z' \in Z(\Omega)$, and $Z \subset Z'$, then $Z' \in \mathcal{F}$.
	\end{enumerate}
\end{defi}
\begin{exa} \label{filter}
	\begin{enumerate}[(a)]
		\item $\mathcal{F} = \{\Omega\}$ is a $z$-filter.
		\item Let $S \subset \Omega$, then 
		\[
		\mathcal{F}(S) = \{F \in Z(\Omega): S \subset F\}
		\]
		is a $z$-filter.
		\item In general, if $\mathfrak{B} \subset 2^{\Omega}$ with the finite intersection property (i.e., if $B_1, B_2 \in \mathfrak{B}$, then $B_1 \cap B_2 \in \mathfrak{B}$), then
		\[
		\mathcal{F}(\mathfrak{B}) = \{F \in Z(\Omega): \exists B \in \mathfrak{B} \text{ such that } B \subset F \}
		\]
		is a $z$-filter.
		\item \label{filterD} Particularly, if $\mathfrak{B} \subset Z(\Omega)$ with the finite intersection property, then there is a ``maximal'' $z$-filter $\mathcal{F} = \mathcal{F}(\mathfrak{B})$ containing $\mathfrak{B}$ (i.e., $\mathfrak{B} \subset \mathcal{F}$), that is, 
		\[
		\mathcal{F} = \{F \in Z(\Omega): \exists B \in \mathfrak{B} \text{ such that } B \subset F \}.
		\]
	\end{enumerate}
\end{exa}

In order to deal with the complex-valued functions case, we need the following definition of antisymmetric $z$-filter, which is analogous to antisymmetric set (see e.g. \cite[Definition 2.5]{Che18f}).
\begin{defi}[antisymmetric $z$-filter]\label{anti}
	Let $v \in F(\Omega, \mathbb{R}_+)$, $A \subset C(\Omega, \mathbb{C})$, and $\mathcal{F}$ a $z$-filter on $\Omega$. We call $\mathcal{F}$ an \emph{$(A, v)$-antisymmetric $z$-filter}, if for each $\phi \in A$ such that
	\[
	\phi(\mathcal{F} \cap \supp v) := \bigcap_{F \in \mathcal{F}} \overline{\phi(F \cap \supp v)} \subset [0, 1],
	\]
	then $\phi(\mathcal{F} \cap \supp v)$ is single\footnote{In particular, here we technically assume that every member of $\mathcal{F}$ meets $\supp v$, i.e., for each $F \in \mathcal{F}$, $F \cap \supp v \neq \emptyset$.}; here, $\overline{\phi(F \cap \supp v)}$ denotes the closure of $\phi(F \cap \supp v)$ in $\mathbb{C}\cup\{\infty\}$. If $\supp v = \Omega$, we also say $\mathcal{F}$ is an \emph{$A$-antisymmetric $z$-filter}.
\end{defi}

\begin{defi}[multiplier]
	For $ W \subset F(\Omega, X) $ and $ \varphi \in F(\Omega, \mathbb{C}) $, we say $ \varphi $ is a \emph{multiplier} of $ W $ if $ \varphi f + (1 - \varphi) g \in W $ for all $ f, g \in W $. A subset $ A \subset F(\Omega, \mathbb{C}) $ is called a \emph{multiplier} of $ W $ if for every $ \varphi \in A $, $ \varphi $ is a multiplier of $ W $.
\end{defi}

\subsection{The distance formula} 
The following is our main result concerning the calculation of $d_{v, p, \Omega}(f, W)$, which generalizes the Machado theorem \cite{Mac77} and \cite[Theorem 2.8]{Che18f}. 

\begin{thm}\label{machado}
	Let $V \subset F(\Omega, \mathbb{R}_+)$, $A \subset C(\Omega, \mathbb{C})$, $W \subset F(\Omega, X)$, and $f \in F(\Omega, X)$. Assume $A$ is a multiplier of $W$. If $f - W \subset FV_{b}(\Omega, X)$, then for each $p \in \mathfrak{A}$, $v \in V$, there is an $(A, v)$-antisymmetric $z$-filter $\mathcal{F}$ such that 
	\[
	\inf_{F \in \mathcal{F}} d_{v, p, F}(f, W) = d_{v, p, \Omega}(f, W).
	\]
	Moreover, the above $z$-filter $\mathcal{F}$ can be chosen as ``maximum'' in the sense that, if there is a $z$-filter $\mathcal{E}$ such that $\mathcal{F} \subset \mathcal{E}$ and
	\[
	\inf_{E \in \mathcal{E}} d_{v, p, E}(f, W) = d_{v, p, \Omega}(f, W),
	\]
	then $\mathcal{E} = \mathcal{F}$.
\end{thm}
\begin{proof}
	The steps of the proof are mainly due to Machado \cite{Mac77} with an ingenious argument due to Brosowski and Deutsch \cite{BD81} (see also \cite{Ran84, Che18f}). 
	
	Let $a = d_{v, p, \Omega}(f, W)$. Note that $a < +\infty$ (as $f - W \subset FV_{b}(\Omega, X)$).
	
	\textbf{Step I}. Consider the collection $\mathbb{P}$ of all $z$-filters $\mathcal{F}$ on $\Omega$ satisfying
	\[
	\inf_{F \in \mathcal{F}} d_{v, p, F}(f, W) = a;
	\]
	particularly, for each $F \in \mathcal{F}$, $d_{v, p, F}(f, W) = a$.
	
	$\mathbb{P}$ is not empty since $\{\Omega\} \in \mathbb{P}$. Partially order $\mathbb{P}$ by inclusion. Then by Hausdorff's maximality theorem\footnote{It seems that unless $A \subset C(\Omega, [0, 1])$ or $A$ is a self-adjoint subalgebra of $C_b(\Omega, \mathbb{C})$, in the general case ($A \subset C(\Omega, \mathbb{C})$), the axiom of choice is needed in the proof; for example, the original proof of the Bishop theorem uses the transfinite induction process (see \cite{Bis61}), the simplified proof of this theorem (see \cite[Theorem 5.7]{Rud91}) due to L. de Branges and I. Glicksberg uses the Krein--Milman theorem, and the original proof of the Machado theorem (see \cite[Theorem 2]{Mac77}) or the simplified one (see e.g. \cite{BD81, Ran84}) uses the Zorn lemma.} (see e.g. \cite[Appendix A]{Rud91}), there is a maximal totally ordered subcollection $\mathbb{L}$ of $\mathbb{P}$. Let $\mathcal{F}_0$ be the union of all the members of $\mathbb{L}$. 
	
	\textbf{Step II}. We show $\mathcal{F}_0 \in \mathbb{L}$.
	$\mathcal{F}_0$ is a $z$-filter, which can be seen easily as follows:
	\begin{enumerate}[(i)]
		\item $\emptyset \notin \mathcal{F}_0$;
		\item for $Z_1, Z_2 \in \mathcal{F}_0$, since $\mathbb{L}$ is totally ordered, there is a $\mathcal{L} \in \mathbb{L}$ such that $Z_1, Z_2 \in \mathcal{L}$, and so $Z_1 \cap Z_2 \in \mathcal{L}$, showing $Z_1 \cap Z_2 \in \mathcal{F}_0$;
		\item if $Z \in \mathcal{F}_0$, then there is $\mathcal{F} \in \mathbb{L}$ so that $Z \in \mathcal{F}$, and so if $Z \subset Z'$, we have $Z' \in \mathcal{F}$ and thus $Z' \in \mathcal{F}_0$.
	\end{enumerate}
	As for each $F \in \mathcal{F}_0$, $d_{v, p, F}(f, W) = a$, we have $\inf_{F \in \mathcal{F}_0} d_{v, p, F}(f, W) = a$. Since $\mathbb{L}$ is maximal totally ordered and for each $\mathcal{F} \in \mathbb{L}$, $\mathcal{F} \subset \mathcal{F}_0$, we must have $\mathcal{F}_0 \in \mathbb{L}$.
	
	\textbf{Step III}. We claim that $\mathcal{F}_0$ is $(A, v)$-antisymmetric. Otherwise, there is $h \in A$ such that
	\[
	h(\mathcal{F}_0 \cap \supp v) := \bigcap_{F \in \mathcal{F}_0} \overline{h(F \cap \supp v)} \subset [0, 1],
	\]
	but $h(\mathcal{F}_0 \cap \supp v)$ is not single. 
	
	Let
	\[
	b = \min h(\mathcal{F}_0 \cap \supp v), ~ c = \max h(\mathcal{F}_0 \cap \supp v).
	\]
	Then $0\leq b < c \leq 1$ and 
	$h(\mathcal{F}_0 \cap \supp v) \subset [b, c]$. 
	Take integers $c_1, c_2, c_3$ such that 
	\[
	0 \leq b<1/c_1<1/c_2<1/c_3<c \leq 1.
	\]
	For $\epsilon > 0$ but small (e.g., $\epsilon < \min\{1/c_1 - b, c-1/c_3\} / 2$), let
	\[
	I^{\inf}_{\epsilon} = [b-\epsilon, 1/c_3] \times [-\epsilon, \epsilon] \subset \mathbb{C}, ~I^{\sup}_{\epsilon} = [1/c_1, c+\epsilon] \times [-\epsilon, \epsilon] \subset \mathbb{C}.
	\]
	Note that $I^{\inf}_{\epsilon} \cup I^{\sup}_{\epsilon} = [b - \epsilon, c+ \epsilon] \times [-\epsilon, \epsilon]$ is a (closed) neighborhood of $[b, c]$ in $\mathbb{C}$. As $h_F := \overline{h(F \cap \supp v)}$ is a compact set of $\mathbb{C}\cup\{\infty\}$, there are finite $G_i \in \mathcal{F}_0$ ($i = 1,2,...,n$) such that $\bigcap h_{G_i} \subset I^{\inf}_{\epsilon} \cup I^{\sup}_{\epsilon}$. Let $F_{0, \epsilon} = \cap G_i \in \mathcal{F}_0$. Then 
	\[
	\overline{h(F_{0, \epsilon} \cap \supp v)} \subset \bigcap h_{G_i} \subset I^{\inf}_{\epsilon} \cup I^{\sup}_{\epsilon},
	\]
	i.e., $F_{0, \epsilon} \cap \supp v \subset h^{-1}(I^{\inf}_{\epsilon} \cup I^{\sup}_{\epsilon})$. Let 
	\[
	O^{\inf}_{\epsilon} = h^{-1}(I^{\inf}_\epsilon), ~O^{\sup}_{\epsilon} = h^{-1}(I^{\sup}_\epsilon).
	\]
	$O^{\inf}_{\epsilon}, O^{\sup}_{\epsilon}$ are zero-sets of $\Omega$.
	
	For $F \in \mathcal{F}_0$, set
	\[
	F^{\inf}_{\epsilon} = F \cap O^{\inf}_{\epsilon}, ~F^{\sup}_{\epsilon} = F\cap O^{\sup}_{\epsilon}.
	\]
	Note that $\emptyset \neq F^{\inf}_{\epsilon} \subsetneq F$, $\emptyset \neq F^{\sup}_{\epsilon} \subsetneq F$ (due to $b < c$ and $\epsilon$ small), and $F^{\inf}_{\epsilon}, F^{\sup}_{\epsilon}$ are zero-sets of $\Omega$. 
	Let
	\[
	\mathfrak{B}^{\inf}_{\epsilon} = \{F^{\inf}_{\epsilon}: F \in \mathcal{F}_0\}, ~\mathfrak{B}^{\sup}_{\epsilon} = \{F^{\sup}_{\epsilon}: F \in \mathcal{F}_0\}.
	\]
	By the finite intersection property of $\mathfrak{B}^{\inf}_{\epsilon}, \mathfrak{B}^{\sup}_{\epsilon}$, there are $z$-filters $\mathcal{F}^{\inf}_{\epsilon}, \mathcal{F}^{\sup}_{\epsilon}$ containing $\mathfrak{B}^{\inf}_{\epsilon}, \mathfrak{B}^{\sup}_{\epsilon}$ (see e.g. \autoref{filter} \eqref{filterD}), respectively.
	Since $F^{\inf}_{\epsilon} \subsetneq F$ and $F^{\inf}_{\epsilon}$ is a zero-set, we see that $\mathcal{F}_0 \subsetneq \mathcal{F}^{\inf}_{\epsilon}$; similarly $\mathcal{F}_0 \subsetneq \mathcal{F}^{\sup}_{\epsilon}$. As $\mathcal{F}_0$ is maximal in $\mathbb{L}$, we get
	\[
	\inf_{F \in \mathcal{F}^{\inf}_{\epsilon}} d_{v, p, F}(f, W) < a, ~\inf_{F \in \mathcal{F}^{\sup}_{\epsilon}} d_{v, p, F}(f, W) < a.
	\]
	
	Fix a small $\epsilon_0 > 0$ (e.g., $\epsilon_0 = \min\{1/c_1 - b, c-1/c_3\} / 4$), and set
	\[
	a_0 = \max\left\{\inf_{F \in \mathcal{F}^{\inf}_{\epsilon_0}} d_{v, p, F}(f, W), \inf_{F \in \mathcal{F}^{\sup}_{\epsilon_0}} d_{v, p, F}(f, W)\right\}.
	\]
	Take a constant $\delta_0 > 0$ such that
	\[
	(1+\delta_0)(a_0 + \delta_0) < a.
	\]
	Take $F_1, F_2 \in \mathcal{F}_0$ such that 
	\[
	d_{v, p, F_1 \cap O^{\inf}_{\epsilon_0}}(f, W) < a_0 + \delta_0, ~d_{v, p, F_2 \cap O^{\sup}_{\epsilon_0}}(f, W) < a_0 + \delta_0.
	\]
	Thus there are $g_{\inf}, g_{\sup} \in W$ such that
	\[
	|f - g_{\inf}|_{v, p, F_1 \cap O^{\inf}_{\epsilon_0}} < a_0 + \delta_0, ~|f - g_{\sup}|_{v, p, F_2 \cap O^{\sup}_{\epsilon_0}} < a_0 + \delta_0.
	\]
	Take
	\[
	C = \max\{|f - g_{\inf}|_{v, p, F_2}, |f - g_{\sup}|_{v, p, F_1}\};
	\]
	note that $C < +\infty$ due to $f - g_{\inf}, f - g_{\sup} \in FV_b(\Omega, X)$ (as $f - W \subset FV_{b}(\Omega, X)$).
	
	Let
	\[
	K_{n}(x) = (1 - x^{n})^{c_2^n}.
	\]
	Note that if $0\leq x \leq 1/c_1$,
	\[
	K_{n}(x) = (1 - x^n)^{c_2^n} \geq 1 - c_2^n x^n \geq 1 - (c_2/c_1)^n;
	\]
	and if $1/c_3 \leq x \leq 1$,
	\[
	K_n(x) = (1-x^n)^{c_2^n} \leq (1+x^n)^{-c_2^n} \leq (c^n_2x^n)^{-1} \leq (c_3/c_2)^n.
	\]
	So we can choose a large integer $N$ such that 
	\[
	\left\{
	\begin{array}{ll}
		1 - K_N(x) < \gamma_0, & x \in [0, 1/c_1],\\
		K_N(x) < \gamma_0, & x \in [1/c_3, 1],
	\end{array}
	\right.
	\]
	where 
	\[
	\gamma_0 = \frac{1}{4}\frac{a - a_0 - \delta_0}{a_0 + \delta_0 + C}.
	\]
	Then we find two neighborhoods $O_1, O_2$ of $[b, 1/c_1], [1/c_3, c]$ in $\mathbb{C}$, respectively, such that
	\[
	\left\{
	\begin{array}{ll}
		|1 - K_N(z)| \leq \gamma_0, & z \in O_1,\\
		|K_N(z)| \leq \gamma_0, & z \in O_2. \\
	\end{array}
	\right.
	\]
	Choose a small $\epsilon > 0$ ($\epsilon < \epsilon_0$) such that for all $z \in I^{\inf}_\epsilon \cup I^{\sup}_\epsilon$,
	\[
	|z| + |1 - z| \leq 1 + \delta_0,
	\]
	and $I^{\inf}_\epsilon \cup I^{\sup}_\epsilon \subset O_1 \cup O_2$, $I^{\inf}_\epsilon \setminus I^{\sup}_\epsilon \subset O_1$, $I^{\sup}_\epsilon \setminus I^{\inf}_\epsilon \subset O_2$. From now on fix such $\epsilon$.
	
	Let 
	\[
	F_* = F_1 \cap F_2 \cap F_{0, \epsilon} \in \mathcal{F}_{0}.
	\]
	Consider
	\[
	g_{N} = h_N \cdot g_{\inf} + (1 - h_N) \cdot g_{\sup}, ~h_N = K_{N}(h),
	\]
	As $h$ is a multiplier of $W$, $g_{N} \in W$\footnote{Here we use the fact: If $ \varphi $ is a multiplier of $ W $, so are $ \varphi^n $, $ 1 - \varphi^n $, and particularly $ (1 - \varphi^n)^m $ where $ n, m \in \mathbb{N}_+ $. In fact (take $ n = 2 $ as an example), for all $ f, g \in W $, since $ \varphi f + (1 - \varphi) g \in W $, we have
		\[
		\varphi^2 f + (1 - \varphi^2) g = \varphi (\varphi f + (1 - \varphi) g) + (1 - \varphi) g \in W.
		\]}. 
	For $t \in F_* \cap \{O^{\inf}_\epsilon \setminus O^{\sup}_\epsilon\}$, $h(t) \in I^{\inf}_\epsilon \setminus I^{\sup}_\epsilon \subset O_1 $, and so
	\begin{align*}
		 & |f - g_N|_{v, p, F_* \cap \{O^{\inf}_\epsilon \setminus O^{\sup}_\epsilon\}} \\
		\leq~ & \sup_{z \in O_1}|K_N(z)||f - g_{\inf}|_{v, p, F_* \cap (O^{\inf}_\epsilon \setminus O^{\sup}_\epsilon)} + \sup_{z \in O_1}|1 - K_N(z)||f - g_{\sup}|_{v, p, F_* \cap (O^{\inf}_\epsilon \setminus O^{\sup}_\epsilon)} \\
		\leq~ & \sup_{z \in O_1}|K_N(z)||f - g_{\inf}|_{v, p, F_1 \cap O^{\inf}_{\epsilon_0}} + \sup_{z \in O_1}|1 - K_N(z)||f - g_{\sup}|_{v, p, F_1} \\
		\leq~ & (1 + \gamma_0)(a_0 + \delta_0) + \gamma_0 C < a.
	\end{align*}
	Similarly,
	\[
	|f - g_N|_{v, p, F_* \cap \{O^{\sup}_\epsilon \setminus O^{\inf}_\epsilon\}} < a.
	\]
	For $t \in F_* \cap \{O^{\inf}_\epsilon \cap O^{\sup}_\epsilon\}$, $h(t) \in I^{\inf}_\epsilon \cap I^{\sup}_\epsilon$, and so
	\begin{align*}
		& v(t)p(f(t) - g_N(t)) \\
		\leq~ & |K_N(h(t))|v(t)p(f(t) - g_{\inf}(t)) + |1 - K_N(h(t))|v(t)p(f(t) - g_{\sup}(t)) \\
		\leq~ & |K_N(h(t))||f- g_{\inf}|_{v, p, F_1 \cap O^{\inf}_{\epsilon_0}} + |1 - K_N(h(t))||f- g_{\sup}|_{v, p, F_2 \cap O^{\sup}_{\epsilon_0}} \\
		\leq~ & (|K_N(h(t))| + |1 - K_N(h(t))|)(a_0 + \delta_0) \\
		\leq~ & (1 + \delta_0) (a_0 + \delta_0) < a,
	\end{align*}
	i.e.,
	\[
	|f - g_N|_{v, p, F_* \cap \{O^{\sup}_\epsilon \cap O^{\inf}_\epsilon\}} < a.
	\]
	We have shown that
	\[
	|f - g_N|_{v, p, F_* \cap \{O^{\sup}_\epsilon \cup O^{\inf}_\epsilon\}} < a.
	\]
	But by the choice of $F_{0, \epsilon}$ (i.e., $F_{0, \epsilon} \cap \supp v \subset h^{-1}(I^{\inf}_{\epsilon} \cup I^{\sup}_{\epsilon}) = O^{\sup}_\epsilon \cup O^{\inf}_\epsilon $), we see that
	\[
	F_* \cap \{O^{\sup}_\epsilon \cup O^{\inf}_\epsilon\} \cap \supp v = F_* \cap \supp v,
	\]
	yielding
	\begin{align*}
		& |f - g_N|_{v, p, F_*} = |f - g_N|_{v, p, F_* \cap \supp v} \\
		=~ & |f - g_N|_{v, p, F_* \cap \{O^{\inf}_\epsilon \cup O^{\sup}_\epsilon\} \cap \supp v} = |f - g_N|_{v, p, F_* \cap \{O^{\inf}_\epsilon \cup O^{\sup}_\epsilon\}} < a.
	\end{align*}
	This is a contradiction to $d_{v, p, F_*}(f, W) = a$.
	The $z$-filter $\mathcal{F}_0$ satisfies the conclusion of the theorem, and the proof is complete.
\end{proof}

Let us consider a situation that in \autoref{machado}, $S := \bigcap_{F \in \mathcal{F}} F \neq \emptyset$, and
\[
d_{v, p, S}(f, W) = d_{v, p, \Omega}(f, W).
\]
\begin{defi}[antisymmetric set]\label{def:anti}
	Let $v \in F(\Omega, \mathbb{R}_+)$, $A \subset C(\Omega, \mathbb{C})$. A subset $ S \subset \supp v $ is called an \emph{$ (A, v) $-antisymmetric set} if $ f \in A $ such that $ f: S \to [0,1] $, then $ f|_{S} $ is a constant. If $\supp v = \Omega$, we also say $S$ is an \emph{$ A $-antisymmetric set}.
\end{defi}
\begin{defi}[bounded condition]\label{def:bc}
	We say a subset $ A $ of $ C(\Omega, \mathbb{C}) $ satisfies the \emph{bounded condition} with respect to $ V $, if every $ \varphi \in A $ is bounded on the support of every $ v \in V $ (see e.g. \cite{Nac65, Pro71, Sum71}). 
\end{defi}
In the classical definition of antisymmetric set (see e.g. \cite{Bis61, Pro71, Pro88}), if
\begin{enumerate}[(a)]
	\item $A$ is a subalgebra and $A \subset C_{b}(\Omega, \mathbb{C})$ (\cite{Bis61}),
	\item or, $A \subset C(\Omega, [0, 1])$ (\cite{Pro88}),
	\item or, $A$ is a subalgebra satisfying the bounded condition with respect to $ \{v\} $ (\cite{Pro71}), 
\end{enumerate}
and if $ A_S = \{ f|_{S}: f\in A \} $ contains no nonconstant real functions (with $S \subset \supp v$ in (c)), then $S$ is called an antisymmetric set. These are special cases of \autoref{def:anti} (see e.g. \cite[Remark 2.6]{Che18f}).

\begin{cor}[See {\cite[Theorem 2.8]{Che18f}}]\label{cor:mac}
	Let $V \subset F(\Omega, \mathbb{R}_+)$ be a family of upper semicontinuous functions, $A \subset C(\Omega, \mathbb{C})$, $W \subset F(\Omega, X)$, and $f \in F(\Omega, X)$. Assume $A$ is a multiplier of $W$. If $f - W \subset CV_{0}(\Omega, X)$, then for each $p \in \mathfrak{A}$, $v \in V$, there is an $(A, v)$-antisymmetric set $S$ such that 
	\[
	d_{v, p, S}(f, W) = d_{v, p, \Omega}(f, W).
	\]
\end{cor}
\begin{proof}
	The direct proof given in \cite[Theorem 2.8]{Che18f} is simple. Here, we give a second proof by \autoref{machado}.
	
	Let $a = d_{v, p, \Omega}(f, W)$. If $a = 0$, then there is nothing to prove. So assume $a > 0$.
	
	Let $\mathcal{F}$ be the $(A, v)$-antisymmetric $z$-filter giving in \autoref{machado}. Let
	\[
	S = \bigcap_{F \in \mathcal{F}} F.
	\]
	We will show that $S \neq \emptyset$ is exactly what we need.
	
	Take $g \in W$. For any $\varepsilon > 0$ ($\varepsilon < a$), $F \in \mathcal{F}$, set
	\[
	\Theta^\varepsilon_F = \{ x \in F: v(x)p(f(x) - g(x)) \geq a - \varepsilon \} \subset \supp v.
	\]
	Then $\Theta^\varepsilon_F \neq \emptyset$. Since $v, p, f - g$ are upper semicontinuous, $\Theta^\varepsilon_F$ is closed, and since $f - g \in CV_{0}(\Omega, X)$, by the definition, $\Theta^\varepsilon_F$ is compact. For any finite $F_i \in \mathcal{F}$ ($i = 1, 2, ..., n$), let $F_* = \bigcap F_i \in \mathcal{F}$, then
	\[
	\bigcap \Theta^\varepsilon_{F_{i}} = \Theta^\varepsilon_{F_*} \neq \emptyset.
	\]
	Thus $\emptyset \neq  \bigcap_{F \in \mathcal{F}} \Theta^\varepsilon_{F} \subset S \cap \supp v$. It follows that for any $g \in W$,
	\[
	\{x \in S : v(x)p(f(x) - g(x)) \geq a - \varepsilon\} \neq \emptyset.
	\]
	This gives that $d_{v, p, S}(f, W) \geq a - \varepsilon$. And so due to $\varepsilon$ arbitrarily small, $d_{v, p, S}(f, W) = a$. 
	
	Next, we show $S$ is an $(A, v)$-antisymmetric set. Without loss of generality, assume $S \subset \supp v$, otherwise take $S \cap \supp v$ instead of $S$.
	If there is $\varphi \in A$ such that $\varphi : S \to [0, 1]$, then we need to show $\varphi|_{S}$ is constant.
	Define 
	\[
	\mathcal{F}(S) = \{F \in Z(\Omega): S \subset F\}.
	\]
	Since $\mathcal{F} \subset \mathcal{F}(S)$ and $d_{v, p, S}(f, W) = a$, we have $\mathcal{F} = \mathcal{F}(S)$ by the ``maximum'' of $\mathcal{F}$.
	
	Take any small neighborhood $O$ of $\overline{\varphi(S)}$, and then find finite rectangles $I_i = [a_i, b_i]\times[c_i, d_i]$ such that $\overline{\varphi(S)} \subset \bigcup I_i \subset O$. As $S \subset \bigcup\varphi^{-1}(I_i)$ and $\bigcup\varphi^{-1}(I_i)$ is a zero-set of $\Omega$, we have
	\[
	\bigcup\varphi^{-1}(I_i) \in \mathcal{F}(S) = \mathcal{F},
	\]
	which yields
	\[
	\bigcap_{F \in \mathcal{F}} \overline{\phi(F \cap \supp v)} = \overline{\varphi(S)} \subset [0, 1].
	\]
	As $\mathcal{F}$ is $(A, v)$-antisymmetric, $\overline{\varphi(S)}$ is single, i.e., $\varphi|_{S}$ is constant. The proof is complete.
\end{proof}

\subsection{Basic facts}

In order to make our results hold for any topological space $\Omega$, we need a result due to Stone and \^{C}ech (see e.g. \cite[Theorem 3.9]{GJ60}).
\begin{thm}[Stone--\^{C}ech]\label{thm:sc}
	For any topological space $ \mathcal{T} $, there exist a completely regular Hausdorff space $ \mathfrak{T} $ and a continuous map $ \pi $ of $ \mathcal{T} $ onto $ \mathfrak{T} $ such that the following hold.
	\begin{enumerate}[(a)]
		\item The map $ g \mapsto g \circ \pi $ is an isomorphism of $ C(\mathfrak{T}, \mathbb{R}) $ onto $ C(\mathcal{T}, \mathbb{R}) $ which carries $C_{b}(\mathfrak{T}, \mathbb{R})$ onto $C_{b}(\mathcal{T}, \mathbb{R})$. Here, $\mathbb{R}$ can be taken as $X$ if $X$ is Hausdorff.
		\item $Z \in Z(\mathcal{T})$ if and only if $\pi(Z) \in Z(\mathfrak{T})$.
	\end{enumerate} 
\end{thm}
In fact, the space is given by $ \mathfrak{T} = \{ [x]_{C(\mathcal{T}, \mathbb{R})}: x \in \mathcal{T} \} $ where 
\[
[x]_{C(\mathcal{T}, \mathbb{R})} = \{ y \in \Omega : g(y) = g(x), \forall g \in C(\mathcal{T}, \mathbb{R}) \}
\]
and the map $ \pi $ is constructed by $ x \mapsto [x]_{C(\mathcal{T}, \mathbb{R})} $. The topology of $ \mathfrak{T} $ is the weakest one such that all $ h: \mathfrak{T} \to \mathbb{R} $ satisfying $ h \circ \pi \in C(\mathcal{T}, \mathbb{R}) $ are continuous. So for finite $ h_{i} \in C(\mathfrak{T}, \mathbb{R}) $ ($i = 1, 2, ..., n$) and $ \varepsilon > 0 $, the set
\[
\{ y: |h_{i}(y) - h_{i}(y_0)| < \varepsilon, i = 1, 2, ..., n \}
\]
is a neighborhood of $ y_0 \in \mathfrak{T} $. In the following, we write 
\begin{enumerate}[$\bullet$]
	\item $ \widehat{\mathcal{T}} = \mathfrak{T} $ if $ \mathfrak{T} $ is constructed through this way,
	\item $\widehat{Z} = \pi(Z)$,
	\item $\widehat{f}$ if $f \in C(\mathcal{T}, \mathbb{R})$ and $\widehat{f} \circ \pi = f$.
\end{enumerate}

Let us give a representation of $\bigcap_{F \in \mathcal{F}} \overline{\phi(F \cap \supp v)}$ in \autoref{anti} through the Stone--\v{C}ech compactification.

Assume $\Omega$ is a completely regular Hausdorff space. Let $\beta \Omega$ be the Stone--\v{C}ech compactification of $\Omega$. For $S \subset \Omega$, $\mathrm{cl}_{\beta} S$ means the closure of $S$ in $\beta \Omega$. For $f \in C(\Omega, X)$ with $\overline{f(\Omega)}$ compact, $f^{\beta}$ denotes the unique Stone extension of $f$ on $\beta \Omega$ such that $f^{\beta}|_S = f$. We will frequently use the following fact: If $Z_1, Z_2 \in Z(\Omega)$, then $\mathrm{cl}_{\beta} (Z_1 \cap Z_2) = \mathrm{cl}_{\beta} Z_1 \cap \mathrm{cl}_{\beta}Z_2$  (see e.g. \cite[6.5 (IV)]{GJ60}). For more details, see \cite[Chapter 6]{GJ60}.

\begin{lem}\label{lem:rep}
	Assume $\Omega$ is a completely regular Hausdorff space and $\mathcal{F}$ is a $z$-filter on $\Omega$. Then for $\phi \in C(\Omega, \mathbb{C})$ (or $\phi \in C(\Omega, X)$ with $\overline{\phi(\Omega)}$ compact),
	\[
	\phi^\beta(\bigcap_{F \in \mathcal{F}}\mathrm{cl}_{\beta}F) = \bigcap_{F \in \mathcal{F}} \phi^\beta(\mathrm{cl}_{\beta}F) = \bigcap_{F \in \mathcal{F}} \overline{\phi(F)}.
	\] 
\end{lem}
\begin{proof}
	First note the following fact: If $A$ is a subset of a Hausdorff topological space with $\overline{A}$ compact and $\varphi$ is a continuous function defined on this space, then $\varphi(\overline{A}) = \overline{\varphi(A)}$. So due to the compactness of $\mathrm{cl}_{\beta}F$ in $\beta\Omega$, 
	\[
	\overline{\phi(F)} = \phi^\beta(\mathrm{cl}_{\beta}F),
	\]
	and (due to the finite intersection property of $\mathcal{F}$)
	\[
	S^\beta := \bigcap_{F \in \mathcal{F}}\mathrm{cl}_{\beta}F \neq \emptyset.
	\]
	Thus,
	\[
	\phi^\beta(S^\beta) \subset \bigcap_{F \in \mathcal{F}} \phi^\beta(\mathrm{cl}_{\beta}F).
	\]
	Let $O$ be any (small) neighborhood of $\phi^\beta(S^\beta)$ in $\mathbb{C} \cup \{\infty\}$ (or $X$). Then $(\phi^\beta)^{-1}(O)$ is an open neighborhood of $S^\beta$. So there are finite $F_i \in \mathcal{F}$ ($i = 1, 2, ..., n$) such that $\bigcap \mathrm{cl}_{\beta}F_i \subset (\phi^\beta)^{-1}(O)$. Since $F_i$ are zero-sets, $\bigcap \mathrm{cl}_{\beta}F_i = \mathrm{cl}_{\beta} \bigcap F_i$, we have $\mathrm{cl}_{\beta} F \subset (\phi^\beta)^{-1}(O)$ where $F := \bigcap F_i \in \mathcal{F}$, i.e., $(\phi^\beta)(\mathrm{cl}_{\beta} F) \subset O$. The proof is complete.
\end{proof}

\begin{rmk}
	\begin{enumerate}[(a)]
		\item From this lemma, we see that when $\Omega$ is a completely regular Hausdorff space, if $\mathcal{F}$ is an $A$-antisymmetric $z$-filter on $\Omega$, then $\bigcap_{F \in \mathcal{F}}\mathrm{cl}_{\beta}F$ is an $A^\beta$-antisymmetric set, where $A^\beta = \{ \phi^{\beta} \in C(\beta \Omega, \mathbb{C} \cup \{\infty\}): \phi \in A \}$.
		\item Let $S \subset \Omega$ be an arbitrary subset. Suppose for each $F \in \mathcal{F}$, $F \cap S \neq \emptyset$. If $\phi^\beta$ is considered as Stone extension of $\phi|_S$ in $\beta S$ and $\mathrm{cl}_{\beta}(F \cap S)$ is the closure taken in $\beta S$, then
		\[
		\phi^\beta(\bigcap_{F \in \mathcal{F}}\mathrm{cl}_{\beta}(F \cap S)) = \bigcap_{F \in \mathcal{F}} \phi^\beta(\mathrm{cl}_{\beta}(F \cap S)) = \bigcap_{F \in \mathcal{F}} \overline{\phi(F \cap S)}.
		\]
		
		\item The lemma also gives the representation of the Stone extension of $f$, i.e., 
		\[
		f^\beta (p) = \bigcap_{F \in \mathcal{F}} \overline{f(F)},
		\]
		where $\mathcal{F}$ is a $z$-ultrafilter (i.e. maximal $z$-filter) with limit $p \in \beta \Omega$.
	\end{enumerate}
\end{rmk}

We say $A \subset C(\Omega, \mathbb{C})$ is self-adjoint, if $\overline{\phi} \in A$ for all $\phi \in A$.

\begin{lem}\label{lem:adjoint}
	Let $v \in F(\Omega, \mathbb{R}_+)$ and $A$ a self-adjoint subalgebra of $C(\Omega, \mathbb{C})$ satisfying the bounded condition with respect to $\{v\}$. Then for each $(A, v)$-antisymmetric $z$-filter $\mathcal{F}$ and $\phi \in A$, the set $\bigcap_{F \in \mathcal{F}} \overline{\phi(F \cap \supp v)}$ is single.
\end{lem}
\begin{proof}
	Without loss of generality, assume $\supp v = \Omega$, otherwise take $\supp v$ instead of $\Omega$ in the following argument. 
	
	(1). First assume $\Omega$ is completely regular and Hausdorff.
	Since $A$ is a self-adjoint subalgebra, the real part and imaginary part of $\phi$ belong to $A$ and so it suffices to consider that $\phi$ is real. Since $\phi$ is bounded, we can assume $-1 \leq \phi \leq 1$. So $\phi^2 \in C(\Omega, [0, 1])$, and by the definition of antisymmetric $z$-filter, $\bigcap_{F \in \mathcal{F}} \overline{\phi^2(F)}$ is single. From \autoref{lem:rep}, we see that $(\phi^\beta)^2|_{F^\beta} = (\phi^2)^\beta|_{F^\beta}$ is constant where $F^\beta = \bigcap_{F \in \mathcal{F}}\mathrm{cl}_{\beta}F$. Let $\varphi = \frac{\phi^3 + \phi^2}{2}$. Note that $\varphi = \frac{\phi + 1}{2} \phi^2 \in A \cap C(\Omega, [0, 1])$. It follows that $\bigcap_{F \in \mathcal{F}} \overline{\varphi(F)}$ is single, i.e., $\frac{(\phi^\beta)^3 + (\phi^\beta)^2}{2}|_{F^\beta}$ is constant. Thus  $\phi^\beta|_{F^\beta}$ is constant, i.e., $\bigcap_{F \in \mathcal{F}} \overline{\phi(F)}$ is single.
	
	(2). When $\Omega$ is not completely regular and Hausdorff, we can use \autoref{thm:sc} to consider the space $\widehat{\Omega}$ and the self-adjoint subalgebra $\widehat{A} = \{ \widehat{\varphi}: \varphi \in A \} \subset C_b(\widehat{\Omega}, \mathbb{C})$. Write $\widehat{\mathcal{F}} = \{\widehat{F} : F \in \mathcal{F}\}$ which is a $z$-filter of $\widehat{\Omega}$. As $\widehat{\varphi}(\widehat{F}) = \varphi(F)$, $\widehat{\mathcal{F}}$ is $\widehat{A}$-antisymmetric. From the first part of the proof, we get $\bigcap_{\widehat{F} \in \widehat{\mathcal{F}}} \overline{\widehat{\phi}(\widehat{F})}$ is single and so is $\bigcap_{F \in \mathcal{F}} \overline{\phi(F)}$. The proof is complete.
\end{proof}

For $A \subset C(\Omega, \mathbb{C})$, we say \emph{$A$ separates disjoint zero-sets} of $\Omega$ if for $Z_1, Z_2 \in Z(\Omega)$ with $Z_1 \cap Z_2 = \emptyset$, there is $\varphi \in A$ such that $\overline{\varphi (Z_1)} \cap \overline{\varphi(Z_2)} = \emptyset$.

\begin{cor}\label{cor:ze} %Assume $\Omega$ is a completely regular Hausdorff space, 
	Assume $\Omega$ is a completely regular Hausdorff space. Let $A \subset C(\Omega, [0, 1])$ or be a self-adjoint subalgebra of $C_{b}(\Omega, \mathbb{C})$. If $A$ separates disjoint zero-sets of $\Omega$, then for any $A$-antisymmetric $z$-filter $\mathcal{F}$, the set $\bigcap_{F \in \mathcal{F}}\mathrm{cl}_{\beta}F$ consists of one point of $\beta \Omega$.
\end{cor}
\begin{proof}
	First note that for each $\phi \in A$, $\phi^\beta(\bigcap_{F \in \mathcal{F}}\mathrm{cl}_{\beta}F)$ is single if $A \subset C(\Omega, [0, 1])$ or $A$ is a self-adjoint subalgebra of $C_{b}(\Omega, \mathbb{C})$ (\autoref{lem:adjoint}). If $p, q \in \bigcap_{F \in \mathcal{F}}\mathrm{cl}_{\beta}F$ with $p \neq q$, then there are two zero-sets $Z_p, Z_q$ of $\Omega$ with $Z_p \cap Z_q \neq \emptyset$ such that $p \in \mathrm{cl}_{\beta}Z_p$, $q \in \mathrm{cl}_{\beta}Z_q$ (see e.g. \cite[6.5(c)]{GJ60}). By the assumption, there is $\varphi \in A$ such that $\overline{\varphi (Z_p)} \cap \overline{\varphi(Z_q)} = \emptyset$; particularly $\varphi^{\beta}(p) \neq \varphi^{\beta}(q)$. This is a contradiction to that $\{\varphi^{\beta}(p), \varphi^{\beta}(q)\} \subset \varphi^\beta(\bigcap_{F \in \mathcal{F}}\mathrm{cl}_{\beta}F)$ is single.
\end{proof}

A well-known result is that for any topological space $\Omega$, disjoint zero-sets are  (completely) separated by $A = C(\Omega, [0, 1])$; e.g., if $\Omega_1 = f^{-1}(0)$ and $\Omega_2 = g^{-1}(0)$ with $\Omega_1  \cap  \Omega_2  = \emptyset $, then $h(x) = \frac{f(x)}{f(x)+g(x)} \in A$ and $h(\Omega_1) = 0$, $h(\Omega_2) = 1$. In some spaces, that $A$ separates disjoint zero-sets of $\Omega$ is the same as that $A$ separates disjoint points of $\Omega$ (i.e., for $x, y \in \Omega$ with $x \neq y$, there is $\phi \in A$ such that $\phi(x) \neq \phi(y)$); for instance, spaces are \emph{almost compact} (see e.g. \citelist{\cite{Hew47} \cite{GJ60}*{6J} \cite{ACIT18}*{Proposition 1.3.9}}).

\section{Applications}

\subsection{An example}
Let us consider a complex-valued uniform approximation problem, which seemingly could not be handled by the results as e.g. \cite{GM93, BM00}.
\begin{exa}\label{exa:rudin}
	Let $I$ be a completely regular and pseudocompact Hausdorff space, $K$ a compact subset of $\mathbb{C}$ which does not separate $\mathbb{C}$ (i.e., $\mathbb{C} \setminus K$ is connected). Let $\varepsilon > 0$. Then for each $f \in C(I \times K, \mathbb{C})$ with $f_t(z) = f(t, z)$ holomorphic in the interior of $K$ (for each $t \in I$), there is $g \in C(I \times K, \mathbb{C})$ such that, $g_t(z) = g(t, z)$ is a polynomial in $z \in K$ (for each $t \in I$), and
	\[
	|f(t, z) - g(t, z)| < \varepsilon
	\]
	for every $(t, z) \in I \times K$.
\end{exa}
\begin{proof}
	Since $I$ is pseudocompact and $K$ is compact, $I \times K$ is pseudocompact and $\beta(I \times K) = \beta I \times K$ by Glicksberg's theorem (see \cite[Theorems 3 and 1]{Gli59} or \autoref{cor:gli}); in particular, $C(I \times K, \mathbb{C}) = C_{b}(I \times K, \mathbb{C})$. Set
	\[
	A = \{ g \in C(I \times K, \mathbb{C}) : g_{t}(z) ~\text{is a polynomial in}~ z \in K \}.
	\]
	We first show that each $A$-antisymmetric $z$-filter $\mathcal{F}$ contains $z$-filters of the form below: 
	\[
	\mathcal{U}_{K} = \{E \times K: E \in \mathcal{U}\}
	\]
	where $\mathcal{U}$ is a $z$-filter on $I$; that is, for each $E \in \mathcal{U}$, $E \times K \in \mathcal{F}$.
	
	Note that $A_0 = C(I, [0, 1])$ separates disjoint zero-sets of $I$. Consider
	\[
	\widetilde{\varphi}(t, z) = \varphi(t) \in A, ~\varphi \in A_0.
	\]
	As $\mathcal{F}$ is $A$-antisymmetric, by \autoref{lem:rep}, $\widetilde{\varphi}^{\beta}(\bigcap _{F \in \mathcal{F}} \mathrm{cl}_{\beta}F)$ is single for each $\widetilde{\varphi}$; and so $\bigcap _{F \in \mathcal{F}} \mathrm{cl}_{\beta}F \subset \{t\} \times K, t \in \beta I$, due to $\bigcap _{F \in \mathcal{F}} \mathrm{cl}_{\beta}F \subset \beta I \times K$. Let $\mathcal{U}$ be the maximal $z$-filter on $I$ with limit $t$. Now for each $E \in \mathcal{U}$, there is $F \in \mathcal{F}$ such that $F \subset E \times K$; particularly $E \times K \in \mathcal{F}$.
	
	Take $f \in C(I \times K, \mathbb{C})$ with $f_t(z) = f(t, z)$ holomorphic in the interior of $K$. We need to show $f \in \overline{A}$, i.e., $d_{I \times K}(f, A) = 0$, and so by \autoref{machado}, only need to show for each $A$-antisymmetric $z$-filter $\mathcal{F}$,
	\[
	\inf_{F \in \mathcal{F}} d_{F}(f, A) = 0.
	\]
	Now it suffices to show 
	\[
	\inf_{E \in \mathcal{U}} d_{E \times K}(f, A) = 0.
	\]
	Consider $f^{\beta} \in C(\beta I \times K, \mathbb{C})$. Since $f$ is bounded and $f^{\beta}$ is obtained through $f$, i.e., for each $t^\beta \in \beta I$, there is $\{t_{\gamma}\} \subset I$ such that $t_{\gamma} \to t^\beta$ and
	\[
	f(t_{\gamma}, z) \to f^{\beta}(t^{\beta}, z)
	\]
	uniformly for $z \in K$ (due to $t \mapsto f^{\beta}(t, \cdot) \in C(K, \mathbb{C})$ continuous and $f^{\beta}|_I = f$ with $I$ dense in $\beta I$), we see that $f^{\beta}(t^\beta, \cdot)$ is holomorphic in the interior of $K$ for each $t^\beta \in \beta I$. Therefore for each $\epsilon > 0$, by Mergelyan's theorem, there is a polynomial $p$ such that
	\[
	|f^{\beta}(t^\beta, z) - p(z)| < \epsilon
	\]
	for every $z \in K$. Let $\mathcal{U}$ be the maximal $z$-filter on $I$ with limit $t^\beta$. There is $E \in \mathcal{U}$ such that
	\[
	|f^{\beta}(t^\beta, z) - f(t, z)| < \epsilon
	\]
	for every $t \in E$, $z \in K$. Thus,
	\[
	|f(t, z) - p(z)| < 2\epsilon
	\]
	for every $t \in E$, $z \in K$, which yields $\inf_{E \in \mathcal{U}} d_{E \times K}(f, A) = 0$. The proof is complete.
\end{proof}
\begin{rmk}
	\begin{enumerate}[(a)]
		\item Let $I$ be any topological space. Then such result also holds for the pseudocompact-open topology (see \autoref{exa:bounded} \eqref{it:ps}). That is, for each $\varepsilon > 0$, each pseudocompact subset $S$ of $I$, and each $f \in C(I \times K, \mathbb{C})$ with $f_t(z) = f(t, z)$ holomorphic in the interior of $K$ (for each $t \in I$), there is $g \in C(I \times K, \mathbb{C})$ such that, $g_t(z) = g(t, z)$ is a polynomial in $z \in K$ (for each $t \in I$), and
		\[
		|f(t, z) - g(t, z)| < \varepsilon
		\]
		for every $(t, z) \in S \times K$.
		\item Suppose $I$ is pseudocompact. Let 
		\[
		\mathcal{H}(K) = \{ f \in C(K, \mathbb{C}) : f ~\text{is holomorphic in the interior of}~ K\},
		\]
		which is a Banach space. $f \in C(I \times K, \mathbb{C})$ if and only if $f \in C(I, C(K, \mathbb{C}))$ (see e.g. \cite[Lemma 1]{Gli59}). And so if $f \in C(I \times K, \mathbb{C})$ with $f_t(z) = f(t, z)$ holomorphic in the interior of $K$ (for each $t \in I$), then $f \in C(I, \mathcal{H}(K))$. Note that $f(I)$ is compact. Now we can use \autoref{thm:SW} to give another proof. 
	\end{enumerate}
\end{rmk}

\subsection{Uniform approximation}\label{sec:uap}

Let us give a generalization of \cite[Section 5]{BM00} to the vector-valued case.
In this subsection, we will give some condition such that $d_{v, p, \Omega}(f, W) = 0$ by $\inf_{F \in \mathcal{F}}d_{v, p, F}(f, W) = 0$ for all $(A, v)$-antisymmetric $z$-filter $\mathcal{F}$ under the situation that $vf$ has relatively compact range (i.e., $\overline{(vf)(\Omega)}$ is compact). % or $\overline{f(\Omega \cap \supp v)}$ is compact
Some examples are:
\begin{enumerate}[(a)]
	\item $f \in CV_0(\Omega, X)$;
	\item $f \in CV_b(\Omega, X)$, and $X$ is one of the following spaces:
	\begin{enumerate}[(1)]
		\item $\mathbb{C}^n$,
		\item $H(O)$ (the holomorphic functions in $O$ where $O$ is an open set of $\mathbb{C}$, see \cite[1.45]{Rud91}),
		\item $C^{\infty}(O, \mathbb{C})$ (the smooth functions in $O$ where $O$ is an open set of $\mathbb{R}^n$, see \cite[1.46]{Rud91}),
		\item $\mathscr{D}_{K}$ (the smooth functions whose compact supports lying in $K$ where $K$ is a compact set of $\mathbb{R}^n$, see \cite[1.46]{Rud91}),
		\item $\mathscr{D}(O)$ (the smooth functions with compact supports in $O$ where $O$ is an open set of $\mathbb{R}^n$, see \cite[6.2 and Theorem 6.5]{Rud91}),
		\item or more generally, $X$ has the Heine--Borel property (i.e., every closed and bounded subset of $X$ is compact, see \cite[page 9]{Rud91});
	\end{enumerate}
	\item $\Omega$ is pseudocompact, $vf \in C(\supp v, X)$ (and so $f \in CV_b(\Omega, X)$), and $X$ is additionally metrizable (i.e., $X$ has a countable local base, see \cite[Theorem 1.24]{Rud91}).
\end{enumerate}

For $\varepsilon > 0$, $x \in X$, $p \in \mathfrak{A}$, we write
\[
B_{p, \varepsilon}(x) := \{y \in X: p(y - x) < \varepsilon\}.
\]

A simple observation in this situation is the following lemma (see also \cite[Theorem 3.1]{BM00}).
\begin{lem}\label{lem:key}
	Let $V \subset F(\Omega, \mathbb{R}_+)$, $f \in F(\Omega, X)$, $W \subset F(\Omega, X)$, $p \in \mathfrak{A}$, $v \in V$. Assume $\overline{(vf)(\Omega)}$ is compact. Then $d_{v, p, \Omega}(f, W) = 0$ if and only if for any $\varepsilon > 0$ and any finite $x_1, x_2, ..., x_n \in X$, there is $g \in W$ such that
	\[
	\text{if}~p(v(x)f(x) - x_i) < \varepsilon, ~\text{then}~p(v(x)g(x) - x_i) < 2\varepsilon, ~ i = 1, 2, ..., n,
	\]
	i.e., $(vf)^{-1}(B_{p, \varepsilon}(x_i)) \subset (vg)^{-1}(B_{p, 2\varepsilon}(x_i))$.
\end{lem}
\begin{proof}
	Necessity. This is obvious (and the compactness of $\overline{(vf)(\Omega)}$ is not needed). Let $\varepsilon > 0$. By $d_{v, p, \Omega}(f, W) = 0$, there is $g \in W$ such that
	\[
	p(v(x)f(x) - v(x)g(x)) = v(x)p(f(x) - g(x)) < \varepsilon.
	\]
	So for any finite $x_1, x_2, ..., x_n \in X$ such that $p(v(x)f(x) - x_i) < \varepsilon$, we have $p(v(x)g(x) - x_i) < 2\varepsilon$.
	
	Sufficiency. Let $\varepsilon > 0$. Since $\overline{(vf)(\Omega)}$ is compact, there are finite $x_1, x_2, ..., x_n \in X$ such that
	\[
	(vf)(\Omega) \subset \bigcup_{i = 1}^{n}B_{p, \varepsilon}(x_i).
	\]
	Let $g$ be the given in the assumption for $\varepsilon > 0, x_1, x_2, ..., x_n$. Now for any $x \in \Omega$, there is $x_i$ such that $p(v(x)f(x) - x_i) < \varepsilon$, and by the assumption we get $p(v(x)g(x) - x_i) < 2\varepsilon$. Thus,
	\[
	v(x)p(f(x) - g(x)) \leq p(v(x)f(x) - x_i) + p(v(x)g(x) - x_i) < 3\varepsilon.
	\]
	This shows that $d_{v, p, \Omega}(f, W) \leq 3\varepsilon$ and so $d_{v, p, \Omega}(f, W) = 0$. The proof is complete.
\end{proof}

We use the notations:

$X^*$: the dual of $X$, that is, $X^*$ consists of all the continuous linear functions from $X$ into $\mathbb{C}$;

$\Re z$: the real part of $z \in \mathbb{C}$.

\begin{defi}\label{def:rho}
	Let $A \subset C(\Omega, \mathbb{C})$, $f \in F(\Omega, X)$, $v \in F(\Omega, \mathbb{R}_+)$. We say ``$A$ \emph{$\beta$-separates disjoint Lebesgue sets of $f$}'' in $\supp v$, if for any $(A, v)$-antisymmetric $z$-filter $\mathcal{F}$, $e^* \in X^*$, $a, b \in \mathbb{R}$ with $a < b$, there is $F \in \mathcal{F}$ such that $ F\cap L_a = \emptyset$ or $F \cap L^b = \emptyset$, where $L_{a}, L^{b}$ are the Lebesgue sets of $e^*f$ in $\supp v$, defined by
	\[
	L_{a} = \{x \in \supp v : \Re e^*(f(x)) \leq a\}, ~L^{b} = \{x \in \supp v : \Re e^*(f(x)) \geq b\}.
	\]
\end{defi}

\begin{lem}\label{lem:sep0}
	Assume $X$ is Hausdorff. Let $A \subset C(\Omega, \mathbb{C})$, $f \in F(\Omega, X)$, $v \in F(\Omega, \mathbb{R}_+)$. Suppose that $\overline{(vf)(\Omega)}$ is compact and $vf$ is continuous on $\supp v$ (resp. $\overline{f(\supp v)}$ is compact and $f$ is continuous on $\supp v$). If $A$ $\beta$-separates disjoint Lebesgue sets of $vf$ (resp. $f$) in $\supp v$, then for any $(A, v)$-antisymmetric $z$-filter $\mathcal{F}$, the set
	\[
	\bigcap_{F \in \mathcal{F}} \overline{(vf)(F \cap \supp v)} \quad \text{(resp.}~  \bigcap_{F \in \mathcal{F}} \overline{f(F \cap \supp v)} \text{)}
	\]
	is single.
\end{lem}
\begin{proof}
	Without loss of generality, assume $\supp v = \Omega$, otherwise take $\supp v$ instead of $\Omega$ in the following argument. Let $\mathcal{F}$ be an $A$-antisymmetric $z$-filter. First, assume $\Omega$ is a completely regular Hausdorff space.  Let $S^\beta = \bigcap_{F \in \mathcal{F}} \mathrm{cl}_{\beta} F$. Since $(vf)^{\beta} (S^{\beta}) = \bigcap_{F \in \mathcal{F}} \overline{(vf)(F)}$ (see e.g. \autoref{lem:rep}), it suffices to show $(vf)^{\beta} (S^{\beta})$ is single. Otherwise, there are $x, y \in S^\beta$ such that $(vf)^{\beta}(x) \neq (vf)^{\beta}(y)$, and so by the Hausdorff property of $X$ there is $e^* \in X^*$ such that
	\[
	a := \Re e^*((vf)^{\beta}(x)) \neq \Re e^*((vf)^{\beta}(y)) := b.
	\]
	Without loss of generality, suppose $a < b$. Let $\varepsilon = (b - a) / 4$. Then there is $F_0 \in \mathcal{F}$ such that $ F_0 \cap L_{a + \varepsilon} = \emptyset$ or $F_0 \cap L^{b - \varepsilon} = \emptyset$ where
	\[
	L_{a + \varepsilon} = \{t \in \Omega : \Re e^*(v(t)f(t)) \leq a + \varepsilon\}, ~L^{b - \varepsilon} = \{t \in \Omega : \Re e^*(v(t)f(t)) \geq b - \varepsilon\}.
	\]
	Note that $x \in \mathrm{cl}_{\beta} L_{a + \varepsilon}$ and $y \in \mathrm{cl}_{\beta} L_{b - \varepsilon}$; and as $L_{a + \varepsilon}$ is a zero-set, $\mathrm{cl}_{\beta} (F \cap L_{a + \varepsilon}) = \mathrm{cl}_{\beta}F \cap \mathrm{cl}_{\beta}L_{a + \varepsilon}$ for all $F \in \mathcal{F}$. Thus 
	\[
	x \in S^\beta \cap \mathrm{cl}_{\beta} L_{a + \varepsilon} = (\bigcap_{F \in \mathcal{F}} \mathrm{cl}_{\beta} F ) \cap \mathrm{cl}_{\beta} L_{a + \varepsilon} = \bigcap_{F \in \mathcal{F}} \mathrm{cl}_{\beta} (F \cap L_{a + \varepsilon});
	\]
	similarly, $y \in S^\beta \cap \mathrm{cl}_{\beta} L^{b - \varepsilon}$. This means that $ F_0 \cap L_{a + \varepsilon} \neq \emptyset$ and $F_0 \cap L^{b - \varepsilon} \neq \emptyset$, which is impossible. 
	
	If $\Omega$ is not completely regular and Hausdorff, using \autoref{thm:sc}, we consider the space $\widehat{\Omega}$, $\widehat{A} = \{ \widehat{\phi} : \phi \in A \} \subset C(\widehat{\Omega}, \mathbb{C})$, and $\widehat{\mathcal{F}} = \{ \widehat{F} : F \in \mathcal{F} \}$ (which is an $\widehat{A}$-antisymmetric $z$-filter). Since $X$ is Hausdorff, there is a unique $\widehat{vf} \in C(\widehat{\Omega}, X)$ such that $vf = \widehat{vf} \circ \pi $. Also note that if $g \in C(\Omega, \mathbb{R})$ then $g(t) \leq a$ if and only if $\widehat{g}([t]) \leq a$, which gives that $\widehat{A}$ $\beta$-separates disjoint Lebesgue sets of $\widehat{vf}$. It follows that $\bigcap_{\widehat{F} \in \widehat{\mathcal{F}}} \overline{(\widehat{vf})(\widehat{F})}$ is single, and so is $\bigcap_{F \in \mathcal{F}} \overline{(vf)(F)}$.
	
	For the statement about $f$, use $f$ instead of $vf$ in the above argument. The proof is complete.
\end{proof}

The following definition for the real-valued functions is widely used in the uniform topology; see e.g. \cite{GM93, BM00}.

\begin{defi}\label{def:sep}
	Let $A \subset C(\Omega, \mathbb{C})$, $f \in F(\Omega, X)$, $S \subset \Omega$. We say ``$A$ \emph{separates disjoint Lebesgue sets of $f$}'' in $S$ if for any $e^* \in X^*$, and the Lebesgue sets of $e^*f$ in $S$,
	\[
	L_{a} = \{x \in S : \Re e^*(f(x)) \leq a\}, ~L^{b} = \{x \in S : \Re e^*(f(x)) \geq b\},
	\]
	where $a, b \in \mathbb{R}$ with $a < b$, there is $\varphi \in A $ such that $\overline{\varphi(L_{a})} \cap \overline{\varphi(L^{b})} = \emptyset$\footnote{Here $\phi(\emptyset) = \emptyset$, $\overline{\emptyset} = \emptyset$ for convenience.}.
\end{defi}

\begin{lem}\label{lem:sep1}
	Let $v \in F(\Omega, \mathbb{R}_+)$ and $f \in F(\Omega, X)$ such that $vf$ (resp. $f$) is continuous on $\supp v$. If $A \subset C(\Omega, [0, 1])$ or $A$ is a self-adjoint subalgebra of $C(\Omega, \mathbb{C})$ satisfying the bounded condition with respect to $ \{v\} $, and $A$ separates disjoint Lebesgue sets of $vf$ (resp. $f$) in $\supp v$, then $A$ $\beta$-separates disjoint Lebesgue sets of $vf$ (resp. $f$) in $\supp v$. 
\end{lem}
\begin{proof}
	Assume $\supp v = \Omega$ (otherwise take $\supp v$ instead of $\Omega$). Without loss of generality, assume $\Omega$ is a completely regular Hausdorff space. Let $\mathcal{F}$ be an $A$-antisymmetric $z$-filter and $S^\beta = \bigcap_{F \in \mathcal{F}} \mathrm{cl}_{\beta} F$. If there exist Lebesgue sets $L_a, L^b$ of $e^*(vf)$ for some $e^* \in X^*$, such that for each $F \in \mathcal{\mathcal{F}}$, $ F\cap L_a \neq \emptyset$ and $F \cap L^b \neq \emptyset$, then we can take
	\[
	x \in \bigcap_{F \in \mathcal{F}} \mathrm{cl}_{\beta} (F \cap L_a) \neq \emptyset, ~y \in \bigcap_{F \in \mathcal{F}} \mathrm{cl}_{\beta} (F \cap L^b) \neq \emptyset.
	\]
	As $L_a$ is a zero-set, $\mathrm{cl}_{\beta} (F \cap L_a) = \mathrm{cl}_{\beta}F \cap \mathrm{cl}_{\beta}L_a$. So $x \in S^\beta \cap \mathrm{cl}_{\beta}L_a$; similarly $y \in S^\beta \cap \mathrm{cl}_{\beta}L^b$. By the assumption that $A$ separates disjoint Lebesgue sets of $vf$, there is $\varphi \in A$ such that $\overline{\varphi(L_{a})} \cap \overline{\varphi(L^{b})} = \emptyset$, i.e., $\varphi^{\beta}(\mathrm{cl}_{\beta} L_a) \cap \varphi^{\beta}(\mathrm{cl}_{\beta} L^b) = \emptyset$. This gives that $\varphi^{\beta}(x) \neq \varphi^{\beta}(y)$. However, since $\mathcal{F}$ is an $A$-antisymmetric $z$-filter, we get $\varphi^\beta(S^{\beta})$ is single (if $A$ is a self-adjoint subalgebra, this follows from \autoref{lem:adjoint}), which is impossible; and thus we finish the proof.
\end{proof}

The following is a generalization of \cite[Theorem 5.4]{BM00} to the complex- and vector-valued functions case; see also \cite[Section 2 Theorems 9 and 11]{Kle75} where the results are stated in $\beta \Omega$.

\begin{thm}\label{thm:ap}
	Let $V \subset F(\Omega, \mathbb{R}_+)$, $A \subset C(\Omega, \mathbb{C})$ and $W \subset FV_{b}(\Omega, X)$. Assume $A$ is a multiplier of $W$. Let $p \in \mathfrak{A}$, $v \in V$, $f \in F(\Omega, X)$. Suppose
	\begin{enumerate}[(i)]
		\item $X$ is Hausdorff;
		\item $\overline{(vf)(\Omega)}$ is compact and $vf$ is continuous on $\supp v$;
		\item $A$ $\beta$-separates disjoint Lebesgue sets of $vf$ in $\supp v$ (see \autoref{def:rho}).
	\end{enumerate}
	Then $d_{v, p, \Omega}(f, W) = 0$ if and only if 
	\begin{enumerate}[  \quad]
		\item  for each $\varepsilon > 0$, $y \in X$, there is $g \in W$ such that $(vf)^{-1}(B_{p, \varepsilon}(y)) \subset (vg)^{-1}(B_{p, 2\varepsilon}(y))$, i.e.,
		\[
		\text{if}~p(v(x)f(x) - y) < \varepsilon, ~\text{then}~p(v(x)g(x) - y) < 2\varepsilon. \label{cond:main} \tag{$\star$}
		\]
	\end{enumerate}
	
\end{thm}
\begin{proof}
	Necessity. This is clear (see also \autoref{lem:key}).
	
	Sufficiency. Without loss of generality, assume $\supp v = \Omega$, otherwise take $\supp v$ instead of $\Omega$. We will show that for each  $A$-antisymmetric $z$-filter $\mathcal{F}$, $\inf_{F \in \mathcal{F}}d_{v, p, F}(f, W) = 0$, and so by \autoref{machado}, $d_{v, p, \Omega}(f, W) = 0$.
	
	Due to $\mathcal{F}$ being $A$-antisymmetric and condition (iii), by \autoref{lem:sep0}, $\bigcap_{F \in \mathcal{F}} \overline{(vf)(F)}$ is single. 
	
	Suppose $\inf_{F \in \mathcal{F}}d_{v, p, F}(f, W) \geq \epsilon_0 > 0$. Let $\bigcap_{F \in \mathcal{F}} \overline{(vf)(F)} = \{ x_0 \}$. Note that each $\overline{(vf)(F)}$ is compact. Then there is $F_0 \in \mathcal{F}$ such that
	\[
	p(v(x)f(x) - x_0) < \epsilon_0 / 8, ~\forall x \in F_0.
	\]
	Now, from \eqref{cond:main}, there is $g \in W$ such that for all $x \in F_0$,
	\[
	p(v(x)g(x) - x_0) < \epsilon_0 / 4,
	\]
	and so for all $x \in F_0$,
	\[
	p(v(x)f(x) - v(x)g(x)) < \epsilon_0 / 8 + \epsilon_0 / 4 < \epsilon_0 / 2,
	\]
	which yields that $d_{v, p, F_0}(f, W) < \epsilon_0 / 2$, a contradiction to $\inf_{F \in \mathcal{F}}d_{v, p, F}(f, W) \geq \epsilon_0$. 
\end{proof}

\begin{rmk}
	If $A$ separates disjoint zero-sets of $\Omega$, then the condition that $X$ is Hausdorff in this theorem is not needed when $\Omega$ is completely regular and Hausdorff, since in this case $\bigcap_{F \in \mathcal{F}} \overline{(vf)(F)} = (vf)^\beta (S^\beta)$ is single due to that $S^\beta$ is single (see \autoref{cor:ze}).
\end{rmk}

\begin{rmk}\label{rmk:Bishop}
	\begin{enumerate}[(a)]
		\item In \cite[Theorem 5.4]{BM00}, Bustamante and Montalvo obtained such result for bounded real-valued functions (i.e. $C_{b}(\Omega, \mathbb{R})$) under an additional assumption that the multiplier $M(W) \subset C(\Omega, [0, 1])$ of $W$ separates disjoint Lebesgue sets of all functions in $W$. This assumption is superfluous and might not be easy to be verified in advance. In \cite[Theorem 4]{Her94}, Hern\'{a}ndez-Mu\~{n}oz also considered such uniform approximation for (real) vector-valued functions in a special case. It seems that the separation condition given in \cite[page 151]{Her94} is restrictive (and so additional assumption that the antisymmetric sets of the multiplier should be single is needed). 
		\item In the complex-valued case (i.e., $A \subset C(\Omega, \mathbb{C})$ and might not be self-adjoint), the separation condition (iii) is needed (instead of that ``$A$ separates disjoint Lebesgue sets'' in the real-valued case); see \autoref{lem:sep1}.
		\item There is no compactness condition on $W$ (and so the Stone extension of the functions $g \in W$ might not have $g^{\beta}(\beta\Omega) \subset X$) in the theorem, which prevents us from deriving such result from the compact case (e.g. \cite{Pro94}); in particular, the method given in \cite{BM00} or \cite{Nel68} cannot be applied to this case.
		\item In general, the compactness condition on $f$ (i.e., condition (ii)) is necessary if we want to use \eqref{cond:main} to deduce $d_{v, p, \Omega}(f, W) = 0$. See \autoref{exa:compact}. The assumption \eqref{cond:main} evidently implies that $d_{v, p, \{x\}}(f, W) = 0$, i.e., $f(x) \in \overline{W(x)}$ for all $x$ with $v(x) > 0$ where $W(x) = \{g(x): g\in W\}$. The latter is widely used in the ``localization problem'' for $\Omega$ being compact or in $CV_0(\Omega, X)$ (see e.g. \cite{Nac65, Pro71, Sum71, Pro77} or \cite[Corollary 3.5]{Che18f}). However, in the uniform topology of $CV_{b}(\Omega, X)$, the condition $d_{v, p, \{x\}}(f, W) = 0$ is generally not enough to give $d_{v, p, \Omega}(f, W) = 0$; intuitively, the asymptotic behavior of $f$ at infinity should be taken into consideration. See \autoref{exa:loc}.
	\end{enumerate}
\end{rmk}

\begin{thm}\label{thm:ap2}
	Let $V \subset F(\Omega, \mathbb{R}_+)$, $A \subset C(\Omega, \mathbb{C})$ and $W \subset FV_{b}(\Omega, X)$. Assume $A$ is a multiplier of $W$. Let $p \in \mathfrak{A}$, $v \in V$, $f \in F(\Omega, X)$. Suppose
	\begin{enumerate}[(i)]
		\item $X$ is Hausdorff;
		\item $\overline{f(\supp v)}$ is compact and $f$ is continuous on $\supp v$;
		\item $A$ $\beta$-separates disjoint Lebesgue sets of $f$ in $\supp v$ (see \autoref{def:rho}).
	\end{enumerate}
	Then $d_{v, p, \Omega}(f, W) = 0$ if and only if for each $\varepsilon > 0$, $y \in X$, there is $g \in W$ such that
	\[
	\text{if}~v(x)p(f(x) - y) < \varepsilon, ~\text{then}~v(x)p(g(x) - y) < 2\varepsilon.
	\]
\end{thm}

\begin{proof}
	Use $f$ instead of $vf$ in the proof of \autoref{thm:ap}.
\end{proof}

The following example shows that the compactness of $vf$ in \autoref{thm:ap} (condition (ii)) usually should be needed.
\begin{exa}\label{exa:compact}
	Let 
	\[
	W = C_0([1, +\infty), C([0, 1], \mathbb{R})) := \{g \in C([1, +\infty), C([0, 1], \mathbb{R})): g(+\infty) = 0\}
	\]
	and $f(x)(y) = y^x$, $x \in [1, +\infty)$, $y \in [0, 1]$. 
	
	$W$ is a closed linear subspace of $C_b([1, +\infty), C([0, 1], \mathbb{R}))$.
	
	Note that $f(+\infty)$ does not exist; particularly $f \notin W $ (in fact $d_{[1, +\infty)}(f, W) = 1$).
	
	Let $A = C([1, +\infty), [0, 1])$. $A$ is a multiplier of $W$ and separates all the zero-sets (i.e. closed sets) of $[1, +\infty)$. In particular, $A$ separates disjoint Lebesgue sets of $f$. We will see below that assumption \eqref{cond:main} in \autoref{thm:ap} holds.
	
	For any $\varepsilon > 0$ ($\varepsilon < 1/4$), $a \in C([0, 1], \mathbb{R})$, consider 
	\[
	\sup_{y \in [0, 1]}|f(x)(y) - a(y)| \leq \varepsilon;
	\]
	particularly, $|f(x)(1) - a(1)| \leq \varepsilon$. As $f(x)(1) = 1$, we see $|1 - a(1)| \leq \varepsilon$ and so $|a(1)| \geq 1/2$. This shows in order that
	\[
	|f(x)(y) - a(y)| = |y^x - a(y)| \leq \varepsilon, ~\forall y \in [0, 1],
	\]
	there must exist $M_0 > 0$ such that $x \leq M_0$ due to $\lim_{x\to +\infty} y^x = 0$, $y \in [0, 1)$. Take $h \in C_{c}([1, +\infty), \mathbb{R})$ with $h|_{[1, M_0]} = 1$. Set
	\[
	g(x)(y) = h(x)g(y).
	\]
	Then $g \in W$ and $|g(x)(y) - a(y)| = 0 < \varepsilon$ for all $y \in [0, 1]$ and $x \leq M_0$.
\end{exa}

The following example gives the difference between $d_{v, p, \{x\}}(f, W) = 0$ and assumption \eqref{cond:main}.

\begin{exa}\label{exa:loc}
	Let $W$ be as in \autoref{exa:compact}. Let $g(x) (y) = 1$ ($\forall y \in [0, 1]$). Then $g \notin W$ (in fact $d_{[1, +\infty)}(g, W) = 1$). But $d_{\{x\}} (g, W)= 0$ for all $x \in [1, +\infty)$, i.e., $g(x) \in W(x)$.
\end{exa}

\subsection{Type of the Stone--Weierstrass theorem}\label{sec:SW}

When $V$ is a Nachbin family, consider the following spaces
\[
CV_{p, 0}(\Omega, X) = \{ f \in CV_{b}(\Omega, X) : \overline{(vf)(\Omega)} \text{ is compact and } vf \in C(\supp v, X) \text{ for each } v \in V \}.
\]
and 
\[
CV_{p, 1}(\Omega, X) = \{ f \in CV_{b}(\Omega, X) : \overline{f(\supp v)} \text{ is compact for each } v \in V \}.
\]
Clearly, if each $v \in V$ is continuous on $\supp v$, then $CV_{p, 0}(\Omega, X)$ equals the classical space 
\[
CV_{p}(\Omega, X) = \{ f \in CV_{b}(\Omega, X) : \overline{(vf)(\Omega)} \text{ is compact for each } v \in V \}.
\]
Moreover, if each $v \in V$ is like $v = \chi_{K}$ for some closed set $K$ (see e.g. \autoref{exa:bounded}), then $CV_{p, 0}$, $CV_{p, 1}$ and $CV_{p}$ are the same.
\begin{thm}\label{thm:w}
	Let $V \subset F(\Omega, \mathbb{R}_+)$ be a Nachbin family. Let $A \subset C(\Omega, \mathbb{C})$ and $W \subset C(\Omega, X)$. Assume
	\begin{enumerate}[(i)]
		\item $\Omega$ is a completely regular Hausdorff space when $X$ is not Hausdorff;
		\item $A \subset C(\Omega, [0, 1])$, or $A$ is a self-adjoint subalgebra satisfying bounded condition with respect to $V$ (see \autoref{def:bc});
		\item $A$ is a multiplier of $W$ (e.g. $W$ is a linear subspace and an $A$-module (i.e. $AW \subset W$)); 
		\item $A$ separates disjoint zero-sets of $\Omega$;
		\item $W$ contains dense constant functions (i.e., there is a dense subset $X_0$ of $X$ such that for each $x \in X_0$, the function $\phi_{x}(t) = x$ ($\forall t \in \Omega$) belongs to W).
	\end{enumerate}
	Then the following hold.
	\begin{enumerate}[(a)]
		\item If $W \subset CV_{p, 0}(\Omega, X)$, then $\overline{W} = CV_{p, 0}(\Omega, X)$ in the topology of $CV_{b}(\Omega, X)$. 
		\item If $W \subset CV_{p, 1}(\Omega, X)$, then $\overline{W} = CV_{p, 1}(\Omega, X)$ in the topology of $CV_{b}(\Omega, X)$. 
	\end{enumerate}
\end{thm}

\begin{proof}
	It suffices to consider (a). Since $A$ separates disjoint zero-sets of $\Omega$, by \autoref{lem:sep1}, $A$ $\beta$-separates disjoint Lebesgue sets of every $f \in C(\Omega, X)$. Since $W$ contains dense constant functions, \eqref{cond:main} is satisfied. And so by \autoref{thm:ap}, for every $f \in CV_{p, 0}(\Omega, X)$, $f \in \overline{W}$. The proof is complete.
\end{proof}

Now we have the following type of the Stone--Weierstrass theorem for $C_b(\Omega, X)$ in the uniform topology when $X$ has the Heine--Borel property or in the pseudocompact-open topology when $X$ is metrizable. Let 
\[
\mathcal{K}C(\Omega, X) = \{ f \in C(\Omega, X) : \overline{f(\Omega)} ~\text{is compact} \}.
\]

\begin{cor}\label{thm:SW}
	Let $A \subset C(\Omega, \mathbb{C})$ and $W \subset C(\Omega, X)$. Assume
	\begin{enumerate}[(i)]
		\item $\Omega$ is a completely regular Hausdorff space when $X$ is not Hausdorff;
		\item $A \subset C(\Omega, [0, 1])$, or $A$ is a self-adjoint subalgebra;
		\item $A$ is a multiplier of $W$ (e.g. $W$ is a linear subspace and an $A$-module (i.e. $AW \subset W$)); 
		\item $A$ separates disjoint zero-sets of $\Omega$;
		\item $W$ contains dense constant functions (i.e., there is a dense subset $X_0$ of $X$ such that for each $x \in X_0$, the function $\phi_{x}(t) = x$ ($\forall t \in \Omega$) belongs to W).
	\end{enumerate}
	Then the following hold.
	\begin{enumerate}[(a)]
		\item If $W \subset \mathcal{K}C(\Omega, X)$ and $A \subset C_{b}(\Omega, \mathbb{C})$, then $\overline{W} = \mathcal{K}C(\Omega, X)$ in the uniform topology. 
		\item If $X$ has the Heine--Borel property, $\overline{W} = C(\Omega, X)$ in the bounded-open topology (\autoref{exa:bounded} \eqref{it:bd}); in particular, if, in addition, $A \subset C_{b}(\Omega, \mathbb{C})$ and $W \subset C_{b}(\Omega, X)$, then $\overline{W} = C_{b}(\Omega, X)$ in the uniform topology.
		\item If $X$ is metrizable, then $\overline{W} = C(\Omega, X)$ in the pseudocompact-open topology.
	\end{enumerate}
\end{cor}
\begin{proof}
	In (b), the Nachbin family $V$ is taken as $V = \{ \chi_{K}: K \text{ is a closed bounded subset of } \Omega \}$, and in (c), $V = \{ \chi_{K}: K \text{ is a closed pseudocompact subset of } \Omega \}$. It follows that $\overline{f(\supp v)}$ is compact where $v \in V, f \in C(\Omega, X)$; furthermore, $A$ satisfies bounded condition with respect to $ V $ (since every continuous function from a bounded (resp. pseudocompact) set into $\mathbb{C}$ is bounded). So it suffices to prove (a), which follows from \autoref{thm:w}. The proof is complete.
\end{proof}

The following example shows that in general, if there exists some function of $C_b(\Omega, X)$ which does not have relatively compact range, then the Stone--Weierstrass theorem does not hold for $C_b(\Omega, X)$ in the uniform topology.

\begin{exa}\label{exa:nonSW}
	Let
	\[
	W_0 = \{ \phi \in C([1, +\infty), C([0, 1], \mathbb{R})) : \overline{\phi([1, +\infty))} ~\text{is compact}  \},
	\]
	and $A =  C([1, +\infty),[0, 1])$. Then $W_0$ is a subalgebra of $C_b([1, +\infty), C([0, 1], \mathbb{R}))$ which is closed in the uniform topology. $A$ is a multiplier of $W_0$, which separates disjoint zero-sets; and $W_0$ contains all the constant functions. But $W_0 \neq $ $C_b([1, +\infty), C([0, 1], \mathbb{R}))$.
\end{exa}

In the following, we consider the problem whether $C(\Omega, X_0)$ is dense in $C(\Omega, X)$ in different topologies where $X_0$ is a dense linear subspace of $X$. Note that \autoref{exa:rudin} is a special case of this problem.

For $\varphi \in C(\Omega, \mathbb{C})$ and $x \in X$, define $\varphi \otimes x$ by
\[
(\varphi \otimes x) (t) = \varphi(t)x.
\]
For $A_0 \subset C(\Omega, \mathbb{C})$ and $X_0 \subset X$, define $A_0 \otimes X_0$ by
\[
A_0 \otimes X_0 = \mathrm{span} \{ \varphi \otimes x : \varphi \in A_0, x \in X_0 \}.
\]

\begin{exa}\label{exa:tensor0}
	Suppose $\Omega$ is a completely regular Hausdorff space when $X$ is not Hausdorff. Let $X_0$ be a dense subset of $X$.
	\begin{enumerate}[(a)]
		\item %Let $\Omega$ be a normal Hausdorff space and $X$ be Hausdorff. 
		$\overline{C_{b}(\Omega, \mathbb{C}) \otimes X_0} = \mathcal{K}C(\Omega, X)$ in the uniform topology.
		\item If $X$ has the Heine--Borel property, then $\overline{C_{b}(\Omega, \mathbb{C}) \otimes X_0} = C_{b}(\Omega, X)$ in the uniform topology, and $\overline{C(\Omega, \mathbb{C}) \otimes X_0} = C(\Omega, X)$ in the bounded-open topology (\autoref{exa:bounded} \eqref{it:bd}).
		\item If $X$ is metrizable, then $\overline{C(\Omega, \mathbb{C}) \otimes X_0} = C(\Omega, X)$ in the pseudocompact-open topology.
	\end{enumerate}
\end{exa}
\begin{proof}
	It suffices to prove (a). Take $A = C(\Omega, [0, 1])$, $W = C_{b}(\Omega, \mathbb{C}) \otimes X_0$. Then $AW \subset W$ and $A$ separates disjoint zero-sets of $\Omega$. $W$ contains dense constant functions. Applying \autoref{thm:SW}, we obtain the result.
\end{proof}

In $CV_0(\Omega, X)$, we also have the similar result.
\begin{thm}\label{thm:litte}
	Let $V \subset C(\Omega, \mathbb{R}_+)$ be a Nachbin family. Let $X_0$ be a dense subset space of $X$. Assume $\Omega$ is a completely regular Hausdorff space. Then $\overline{CV_0(\Omega, \mathbb{C}) \otimes X_0} = CV_0(\Omega, X)$ in the topology of $CV_0(\Omega, X)$. 
\end{thm}
\begin{proof}
	Let $A = C(\Omega, [0, 1])$ and $W = CV_0(\Omega, \mathbb{C}) \otimes X_0$. Then $A W \subset W$ and so $A$ is a multiplier of $W$. Since $\Omega$ is completely regular and Hausdorff, $A$ separates disjoint points of $\Omega$, which means that all the $A$-antisymmetric sets are $\{s\}, s \in \Omega$. Let $f \in CV_0(\Omega, X)$. For each $s \in \Omega$, $g(t) = f(s)$ ($t \in \Omega$) belongs to $\overline{W}$. Therefore, from the Localization Theorem (see \cite[Corollary 3.5]{Che18f}, a corollary of \autoref{cor:mac}), $f \in \overline{W}$. The proof is complete.
\end{proof}

\begin{exa}
	Suppose $\Omega$ is a completely regular Hausdorff space when $X$ is not Hausdorff. Let $X_0$ be a dense subset of $X$.
	\begin{enumerate}[(a)]
		\item $\overline{C_0(\Omega, \mathbb{C}) \otimes X_0} = C_0(\Omega, X)$ in the uniform topology.
		\item $\overline{C_b(\Omega, \mathbb{C}) \otimes X_0} = C_b(\Omega, X)$ in the strict topology.
		\item $\overline{C(\Omega, \mathbb{C}) \otimes X_0} = C(\Omega, X)$ in the compact-open topology.
	\end{enumerate}
\end{exa}
\begin{proof}
	If $\Omega$ is a completely regular Hausdorff space, then the result follows from \autoref{thm:litte} by taking different Nachbin families (see \autoref{exa:c0}). For (a) (c), if $\Omega$ is any topological space but $X$ is Hausdorff, then we can use \autoref{thm:sc} to consider the space $\widehat{\Omega}$, which is completely regular and Hausdorff; we give the details below.
	
	(a). Let $f \in C_0(\Omega, X)$ and $\Omega_0 = f^{-1}(X\setminus\{0\})$. Applying \autoref{thm:sc} to $\Omega_0$, we get the space $\widehat{\Omega}_0$. Since $X$ is Hausdorff, there is a unique $\widehat{f} \in C(\widehat{\Omega}_0, X)$ such that $\widehat{f} \circ \pi = f$.
	As shown in the proof of \cite[Corollary 3.15 (a)]{Che18f}, $\Omega_0$ and $\widehat{\Omega}_0$ are locally compact, $\widehat{f} \in C_0(\widehat{\Omega}_0, X)$, and $\pi:\Omega_0 \to \widehat{\Omega}_0$ is proper. Let $\epsilon > 0$. Write $C_c(\widehat{\Omega}_0, \mathbb{C}) = \{f \in C(\widehat{\Omega}_0, \mathbb{C}) : \supp f \text{ is compact}\}$. Then $C_c(\widehat{\Omega}_0, \mathbb{C}) \otimes X_0$ is dense in $C_0(\widehat{\Omega}_0, X)$ (see e.g. \cite[Example 3.10]{Che18f}). So we have $\widehat{g} \in C_c(\widehat{\Omega}_0, \mathbb{C}) \otimes X_0$ such that $|\widehat{f} - \widehat{g}|_{\widehat{\Omega}_0} < \epsilon$. Set $g = \widehat{g} \circ \pi \in C_c(\Omega_0, X) \subset C_0(\Omega, X)$. Then we obtain 
	\[
	|f - g|_{\Omega} = |f - g|_{\Omega_0} = |\widehat{f} - \widehat{g}|_{\widehat{\Omega}_0} < \epsilon.
	\]
	
	(c). Let $f \in C(\Omega, X)$. Since $X$ is Hausdorff, there is a unique $\widehat{f} \in C(\widehat{\Omega}, X)$ such that $\widehat{f} \circ \pi = f$. Let $K$ be a closed compact subset of $\Omega$. Then $\widehat{K} = \pi(K)$ is a compact subset of $\widehat{\Omega}$. Let $\epsilon > 0$ and $p \in \mathfrak{A}$. Then there is $\widehat{g} \in C(\widehat{\Omega}, \mathbb{C}) \otimes X_0$ such that $|\widehat{f} - \widehat{g}|_{p, \widehat{K}} < \epsilon$. Set $g = \widehat{g} \circ \pi$. Then $|f - g|_{p, K} = |\widehat{f} - \widehat{g}|_{p, \widehat{K}} < \epsilon$. The proof is complete.
\end{proof}

If $\Omega$ and $X$ have nicer topologies, we have the following result.
\begin{thm}
	Let $X_0$ be a dense linear subspace of $X$ and $f \in C(\Omega, X)$.
	If $\Omega$ is paracompact and Hausdorff or $X$ is separable, then for each $\epsilon > 0$, $p \in \mathfrak{A}$ and $f \in C(\Omega, X)$, there is $g \in C(\Omega, X_0)$ such that $|f - g|_{p, \Omega} < \epsilon$.
\end{thm}
\begin{proof}
	(1). First we assume $X$ is separable. Let $\{x_i: i = 1, 2, 3, ...\}$ be a countable dense set in $X$. Set $B_{i}(\epsilon) = \{ x : p(x - x_i) < \epsilon / 2\}$, which is a cozero-set of $X$. Let $U_{i} = f^{-1}(B_{i}(\epsilon))$. Then $\{U_{i}\}$ is countable cover of $\Omega$ by cozero-sets. From \cite[Lemma 2.1]{BH74}, there is a countable locally finite partition of unity $\{h_i\}$ on $X$ with $\supp h_i \subset U_i$. Take $y_i \in U_i \cap X_0$. Set $g(t) = \sum_{i} h_i(t)y_i$. Then $g\in C(\Omega, X_0)$ and 
	\[
	p(f(t) - g(t)) \leq \sum_{i} h_i(t)p(f(t) - y_i) < \epsilon.
	\]
	
	(2). Now we suppose that $\Omega$ is paracompact and Hausdorff. The proof is essentially the same as (1). Set $B_{p, \epsilon} (x) = \{y : p(y - x) < \epsilon\}$ and $U_{x} = f^{-1}(B_{p, \epsilon} (x))$. Then $\{U_{x} : x \in X\}$ is an open cover of $\Omega$. Since $\Omega$ is paracompact, there is a locally finite partition of unity $\{h_x\}$ on $X$ with $\supp h_x \subset U_x$ (see e.g. \cite[Chapter VIII Theorem 4.2]{Dug66}). Take $y_x \in U_x \cap X_0$. Set $g(t) = \sum_{x} h_x(t)y_x$. Then $g\in C(\Omega, X_0)$ and 
	\[
	p(f(t) - g(t)) \leq \sum_{x} h_x(t)p(f(t) - y_x) < 2\epsilon.
	\]
	The proof is complete.
\end{proof}

In the following result, using the classical Stone--Weierstrass theorem, Hager \cite[Proposition 5]{Hag66} proved first that $C(\beta I \times \beta J, \mathbb{R}) = \overline{C_{b}(I, \mathbb{R}) \otimes C_{b}(J, \mathbb{R})} = \overline{C(\beta I, \mathbb{R}) \otimes C(\beta J, \mathbb{R})}$, and then by Glicksberg's theorem, obtaining (a) $\Rightarrow$ (b). Another proof of (a) $\Rightarrow$ (b) by involving Glicksberg's theorem was given in \cite{Tam61}. Let us give a direct proof of (a) $\Rightarrow$ (b) by using our vector-valued version of the Stone--Weierstrass theorem in the uniform topology, but also involving the fundamental Glicksberg's lemma \cite[Lemma 1]{Gli59}.

\begin{cor}[{\cite{Gli59, Tam61}}]\label{cor:gli}
	Let $I, J$ be completely regular, Hausdorff, and infinite. Then the following statements are equivalent. 
	\begin{enumerate}[(a)]
		\item $I \times J$ is pseudocompact.
		\item $C_{b}(I \times J, \mathbb{R}) = \overline{C_{b}(I, \mathbb{R}) \otimes C_{b}(J,  \mathbb{R})}$ in the uniform topology.
		\item $\beta (I \times J) = \beta I \times \beta J$.
	\end{enumerate}
\end{cor}

\begin{proof}
	(a) $\Rightarrow$ (b). 
	By \cite[Lemma 1]{Gli59}, $C(I \times J, \mathbb{R}) = C(I, C(J, \mathbb{R}))$; and since $I$ is pseudocompact and $C(J, \mathbb{R})$ is metrizable, from \autoref{exa:tensor0} (c), we see $C(I, X) = \overline{C(I, \mathbb{R}) \otimes X}$, where $X = C(J, \mathbb{R})$. 
	
	(b) $\Rightarrow$ (c) (see also \cite{Tam61, Hag66}). From the uniqueness of the Stone--\v{C}ech compactification, it suffices to show that for each $f \in C_{b}(I \times J, \mathbb{R})$, there is $\tilde{f} \in C(\beta I \times \beta J, \mathbb{R})$ such that $\tilde{f}|_{I \times J} = f$. Since $C_{b}(I \times J, \mathbb{R}) = \overline{C_{b}(I,  \mathbb{R}) \otimes C_{b}(J, \mathbb{R})}$, there are $g_n \in C_{b}(I, \mathbb{R}) \otimes C_{b}(J, \mathbb{R})$ ($n = 1, 2, 3, ...$) such that $|f - \sum_{i=1}^{n}g_i| < 2^{-n}$. Note that $\sup_{t\in I \times J}|g_n(t)| \leq 3 \cdot 2^{-n-1}$ ($n \geq 2$), $f = \sum_{n=1}^{+\infty} g_n$, and each $g_n$ can be extended continuously to a function $\tilde{g}_n$ on $\beta I \times \beta J$ with $\max_{t\in \beta I \times \beta J}|\tilde{g}_n(t)| = \max_{t\in I \times J}|g_n(t)|$. Define $\tilde{f} = \sum_{n=1}^{+\infty} \tilde{g}_n$ as desired.
	
	(c) $\Rightarrow$ (a) was proved in \cite{Gli59}. 
\end{proof}

\subsection{Splitting of the tensor product type}\label{sec:tensor}

In this subsection, we consider the splitting of $C(I \times J, X \otimes Y)$ as the closure of $C(I, X) \otimes C(J, Y)$ in different senses. See e.g. \citelist{\cite{Tre67}*{Chapters 39--40 (and 44--46)}  \cite{Kle75}*{Section 4}  \cite{Bie81}  \cite{Tim03}}  for some related results. 

Let $X, Y$ be locally convex topological vector (Hausdorff) spaces (over the field $\mathbb{C}$). Let $X \otimes Y$ be the tensor product of $X$ and $Y$. Give a topology on $X \otimes Y$, making $X \otimes Y$ a locally convex topological vector (Hausdorff\footnote{Here, we mean if $X, Y$ are Hausdorff then $X \otimes Y$ should be Hausdorff.}) space and the canonical bilinear map $\otimes: X \times Y \to X \otimes Y, (x, y) \mapsto x \otimes y$ continuous; for example, the $\varepsilon$-topology or the $\pi$-topology on $X \otimes Y$ (see e.g. \cite[Chapter 43]{Tre67} for more details).

For $f \in F(\Omega, X), g \in F(\Omega, Y)$, define $f \otimes g \in F(\Omega, X \otimes Y)$ by
\[
(f \otimes g)(t) = f(t) \otimes g(t), ~ t \in \Omega.
\]
This gives the tensor product of $F(\Omega, X) \otimes F(\Omega, Y)$; particularly, if $W \subset F(\Omega, X)$ and $V \subset F(\Omega, Y)$, then 
\[
W \otimes V = \mathrm{span} \{ f \otimes g : f \in W, g \in V \}.
\]
Since $\otimes: X \times Y \to X \otimes Y$ is continuous, $C_{b}(\Omega, X) \otimes C_{b}(\Omega, Y) \subset C_{b}(\Omega, X \otimes Y)$; moreover, the canonical bilinear map
\[
C_{b}(\Omega, X) \times C_{b}(\Omega, Y) \to C_{b}(\Omega, X) \otimes C_{b}(\Omega, Y), ~(f, g) \mapsto f \otimes g,
\]
is continuous in the uniform topology of $C_{b}(\Omega, X \otimes Y)$. It follows that if $W \subset C_{b}(\Omega, X)$ and $V \subset C_{b}(\Omega, Y)$, then
\[
\overline{W \otimes V} = \overline{\overline{W} \otimes \overline{V}}
\]
in the uniform topology where $\overline{W}$ (resp. $\overline{V}$) is taken in $C_{b}(\Omega, X)$ (resp. $C_{b}(\Omega, Y)$).

\begin{lem}\label{lem:tensor}
	Let $\Omega$ be a completely regular Hausdorff space. Let $X, Y$ be locally convex topological vector spaces. Then,
	\begin{enumerate}[(a)]
		\item $\mathcal{K}C(\Omega, X \otimes Y) = \overline{\mathcal{K}C(\Omega, X) \otimes \mathcal{K}C(\Omega, Y)}$ in the uniform topology;
		\item if $X, Y$ have the Heine--Borel property, then $C_{b}(\Omega, X \otimes Y) = \overline{C_{b}(\Omega, X) \otimes C_{b}(\Omega, Y)}$ in the uniform topology, and $C(\Omega, X \otimes Y) = \overline{C(\Omega, X) \otimes C(\Omega, Y)}$ in the bounded-open topology (\autoref{exa:bounded} \eqref{it:bd}).
		\item if $X, Y$ are metrizable, then $C(\Omega, X \otimes Y) = \overline{C(\Omega, X) \otimes C(\Omega, Y)}$ in the pseudocompact-open topology.
	\end{enumerate}
\end{lem}

\begin{proof}
	Note that if $X, Y$ have the Heine--Borel property, so is $X \otimes Y$ due to that $\otimes: X \times Y \to X \otimes Y$ is continuous; also if $X, Y$ are metrizable, so is $X \otimes Y$.
	In (b), $C_{b}(\Omega, X \otimes Y) = \mathcal{K}C(\Omega, X \otimes Y)$; and in (c), $C(\Omega, X \otimes Y) = \mathcal{K}C(\Omega, X \otimes Y)$ if $\Omega$ is pseudocompact. So it suffices to prove (a). 
	
	Let $A = C(\Omega, [0, 1])$. Then it separates disjoint zero-sets of $\Omega$. Let $W = \mathcal{K}C(\Omega, X) \otimes \mathcal{K}C(\Omega, Y)$. $W$ contains all the constant functions, $A W \subset W$. So by \autoref{thm:SW}, $\overline{W} = \mathcal{K}C(\Omega, X \otimes Y)$ in the uniform topology.
\end{proof}
	
\begin{thm}
	Let $I, J$ be completely regular Hausdorff spaces, and $X, Y$ locally convex topological vector spaces. Assume $I \times J$ is pseudocompact. Then,
	\begin{enumerate}[(a)]
		\item $\mathcal{K}C(I \times J, X \otimes Y) = \overline{\mathcal{K}C(I, X) \otimes \mathcal{K}C(J, Y)}$ in the uniform topology;
		\item if $X, Y$ have the Heine--Borel property or are metrizable, then $C(I \times J, X \otimes Y) = \overline{C(I, X) \otimes C(J, Y)}$ in the uniform topology.
	\end{enumerate}
\end{thm}

\begin{proof}
	Note that in (b), $C(I \times J, X \otimes Y) = \mathcal{K}C(I \times J, X \otimes Y)$, so it suffices to prove (a). By \autoref{lem:tensor}, 
	\[
	\mathcal{K}C(I \times J, X \otimes Y) = \overline{\mathcal{K}C(I \times J, X) \otimes \mathcal{K}C(I \times J, Y)}.
	\]
	By \autoref{exa:tensor0}, 
	\[
	\mathcal{K}C(I \times J, X) = \overline{C(I \times J, \mathbb{C}) \otimes X},
	\]
	and by \autoref{cor:gli}, 
	\[
	C(I \times J, \mathbb{C}) = \overline{C(J, \mathbb{C}) \otimes C(I, \mathbb{C})},
	\]
	which yields 
	\[
	\mathcal{K}C(I \times J, X) = \overline{(C(J, \mathbb{C}) \otimes C(I, \mathbb{C})) \otimes X} = \overline{C(J, \mathbb{C}) \otimes C(I, X)};
	\]
	here, note that $(C(J, \mathbb{C}) \otimes C(I, \mathbb{C})) \otimes X$ is canonically isomorphic onto $C(J, \mathbb{C}) \otimes (C(I, \mathbb{C}) \otimes X)$ (for the topological vector space structure). Similarly, $\mathcal{K}C(I \times J, Y) = \overline{C(I, \mathbb{C}) \otimes C(J, Y)}$. Also, note that $(C(J, \mathbb{C}) \otimes C(I, X)) \otimes (C(I, \mathbb{C}) \otimes C(J, Y))$ is canonically isomorphic onto $C(I, X) \otimes C(J, Y)$ (for the topological vector space structure). It follows that $\mathcal{K}C(I \times J, X \otimes Y) = \overline{\mathcal{K}C(I, X) \otimes \mathcal{K}C(J, Y)}$ and the proof is complete.
\end{proof}

For the sake of completeness, we continue to consider the splitting of tensor products in $CV_0(I \times J, X \otimes Y)$.
For $W \subset F(\Omega, X)$, write
\[
W(s) = \{ f(s): f \in W \}.
\]

\begin{thm}\label{thm:little}
	Let $I, J$ be completely regular Hausdorff spaces, and $X, Y$ locally convex topological vector spaces. Let $V$ (resp. $W$, resp. $U$) be a Nachbin family on $I$ (resp. $J$, resp. $I \times J$). Suppose 
	\begin{enumerate}[(i)]
		\item $\overline{CV_0(I, X)(t) \otimes CW_0(J, Y)(s)} = \overline{CU_{0}(I \times J, X \otimes Y)(t, s)}$ ($\forall (t, s) \in I \times J$), and 
		\item $CV_{0}(I, X) \otimes CW_{0}(J, Y) \subset CU_{0}(I \times J, X \otimes Y)$.
	\end{enumerate}
	Then $CU_{0}(I \times J, X \otimes Y) = \overline{CV_{0}(I, X) \otimes CW_{0}(J, Y)}$ in the topology of $CU_{0}(I \times J, X \otimes Y)$.
\end{thm}

\begin{proof}
	Take $A_{I} = C(I, [0, 1])$ and $A_{J} = C(J, [0, 1])$. Let $A = A_{I}A_{J} := \{ hk : h \in A_{I}, k \in A_{J} \} $, $\mathcal{W} = CV_{0}(I, X) \otimes CW_{0}(J, Y)$. Then $A \mathcal{W} \subset \mathcal{W}$; particularly, $A$ is a multiplier of $\mathcal{W}$. Since $A_I$ (resp. $A_J$) separates $I$ (resp. $J$) (i.e., for any $t_1, t_2 \in I$ with $t_1 \neq t_2$, there is $\varphi \in A_{I}$ such that $\varphi(t_1) \neq \varphi(t_1)$). It follows that every $A$-antisymmetric set reduces to a point. Take $f \in CU_{0}(I \times J, X \otimes Y)$. Since $\overline{CV_0(I, X)(t) \otimes CW_0(J, Y)(s)} = \overline{CU_{0}(I \times J, X \otimes Y)(t, s)}$, we have $f(t, s) \in \overline{\mathcal{W}(t, s)}$ ($\forall (t, s) \in I \times J$). Now, from the Localization Theorem (see \cite[Corollary 3.5]{Che18f}, a corollary of \autoref{cor:mac}), we see $f \in \overline{\mathcal{W}}$, finishing the proof.
\end{proof}

\begin{cor}
	Let $I, J$ be completely regular Hausdorff spaces, and $X, Y$ locally convex topological vector spaces.
	Then,
	\begin{enumerate}[(a)]
		\item if $I \times J$ is locally compact, then $C_{0}(I \times J, X \otimes Y) = \overline{C_{0}(I, X) \otimes C_{0}(J, Y)}$ in the uniform topology;
		\item $C_{b}(I \times J, X \otimes Y) = \overline{C_{b}(I, X) \otimes C_{b}(J, Y)}$ in the strict topology;
		\item $C(I \times J, X \otimes Y) = \overline{C(I, X) \otimes C(J, Y)}$ in the compact-open (and point-open) topology.
	\end{enumerate}
\end{cor}

\begin{proof}
	Since $\otimes: X \times Y \to X \otimes Y$ is continuous, we see $C_{0}(I, X) \otimes C_{0}(J, Y) \subset C_{0}(I \times J, X \otimes Y)$, $C_{b}(I, X) \otimes C_{b}(J, Y) \subset C_{b}(I \times J, X \otimes Y)$ and $C(I, X) \otimes C(J, Y) \subset C(I \times J, X \otimes Y)$. Also note that in (a), since $I \times J$ is locally compact, $C_{0}(I \times J, X \otimes Y) (t, s) = X \otimes Y = C_{0}(I, X) (t) \otimes C_{0}(J, Y) (s)$. Now the statements follow from \autoref{thm:little}.
\end{proof}

\subsection{Type of the Tietze--Dugundji extension theorem} \label{sec:extension}

For $W \subset F(\Omega, K)$, write $W|_{S} = \{f|_{S}: f \in W\}$. The following result is re-proved by many authors especially under the situation that $S$ is a compact subset of a completely regular space or $S$ is a closed subset of a normal space (e.g. \cite[Theorem 3.1]{Pro77} and \cite[Theorem 1]{Tim02}); note that in this case, $C(\Omega, [0, 1])|_S = C(S, [0, 1])$ and $C(S, [0, 1])$ separates disjoint Lebesgue sets of $S$.
\begin{lem}\label{lem:pre}
	Let $f \in C(S, X)$ where $S \subset \Omega$. Assume that $A = C(\Omega, [0, 1])$ separates disjoint Lebesgue sets of $f$ in $S$ (see \autoref{def:sep}) and $\overline{f(S)}$ is compact. Let 
	\[
	W = \left\{\sum_{i = 1}^{n} h_{i}(t)x_i: x_i \in f(S), h_{i} \in C(\Omega, [0, 1]), \sum_{i = 1}^{n} h_{i} = 1, n = 1, 2, 3,... \right\}.
	\]
	Then $f \in \overline{W}$ in $C_{b}(S, X)$.
\end{lem}

\begin{proof} 
	Note that $A = C(\Omega, [0, 1])$ is a multiplier of $W$. Let $p$ be a seminorm of $X$. In order to apply \autoref{thm:ap} to $f \in C_{b}(S, X)$, $W|_{S}$ and $A|_{S}$, obtaining $d_{p, S}(f, W) = 0$, we need to show that assumption \eqref{cond:main} holds.  For each $\epsilon > 0$, $y \in X$, choose $y_0 \in f(S)$ such that $p(y_0 - y) < \epsilon$ (otherwise for such $y$, \eqref{cond:main} holds). Letting $g(t) = y_0$, we see $g \in W$. 
	
	Now for each open set $V$ of $W$ in $C_{b}(S, X)$, we can find a finite number of seminorms $p_{i}$ of $X$, such that $W + V_0 \subset V$ where $ V_0 := \{ f \in C_{b}(S, X) : |f|_{p, S} < \epsilon \} $ is a neighborhood of $0$ in $C_{b}(S, X)$ for some small $\epsilon > 0$ and $p = \max\{p_{i}\}$. Since $d_{p, S}(f, W) = 0$, there is $g \in W$ such that $|f - g|_{p, S} < \epsilon$, that is, $f \in W + V_0$.
	The proof is complete.
\end{proof}

A subset $S$ of $\Omega$ is \emph{$C^*$-embedded} in $\Omega$ if every function in $C_b(S, \mathbb{R})$ can be extended to a function in $C_b(\Omega, \mathbb{R})$. The Gillman--Jerison's version of the Urysohn extension theorem (\cite[Theorem 1.17]{GJ60}) says that $S$ is $C^*$-embedded in $\Omega$ if and only if any two completely separated sets in $S$ are completely separated in $\Omega$ (i.e., for $S_1, S_2 \subset S$, if there is $\phi \in C(S, [0, 1])$ such that $\phi(S_1) = 0, \phi(S_2) = 1$, then there is $\varphi \in C(\Omega, [0, 1])$ such that $\varphi(S_1) = 0, \varphi(S_2) = 1$). That is, $S$ is $C^*$-embedded in $\Omega$ if and only if $A = C(\Omega, [0, 1])$ separates disjoint zero-sets of $S$; see also \cite[Section 3]{BH74} for more details. Typical $C^*$-embedded subsets are: compact subsets of completely regular spaces, and closed subsets of normal spaces. Let us give a vector-valued version of this theorem due to Gutev, Ohta and Yamazaki (in the set-valued case).

\begin{thm}[{See \cite[Theorem 1.3]{GOY06}}]\label{thm:gj}
	 Let $X$ be a Fr\'{e}chet space. $S$ is $C^*$-embedded in $\Omega$ if and only if for every $f \in C(S, X)$ with $f(S)$ relatively compact (i.e., equivalently, totally bounded in Fr\'{e}chet spaces), there exists $g \in C(\Omega, X)$ such that $g|_{S} = f$ and $g(\Omega)$ belongs to the convex hull of $f(S)$.
\end{thm}

\begin{proof}
	Necessity. For $S_1, S_2 \in Z(S)$, if there is $\phi \in C(S, [0, 1])$ such that $\phi(S_1) = 0, \phi(S_2) = 1$, then $\phi \otimes x \in C(S, X)$ where $x \in X$, and $(\phi \otimes x)(S)$ is relatively compact, and so there is $g \in C(\Omega, X)$ such that $g|_{S} = \phi \otimes x$. Then $\varphi = x^*g$, where $x^* \in X^*$ with $x^*(x) = 1$, satisfies $\varphi(S_1) = 0, \varphi(S_2) = 1$.
	
	Sufficiency. Since $S$ is $C^*$-embedded in $\Omega$, $A = C(\Omega, [0, 1])$ separates disjoint zero-sets of $S$, and in particular, separates disjoint Lebesgue sets of every $f \in C(S, X)$ in $S$. Now the conclusion follows from \autoref{lem:pre} and \cite[Lemma 1]{Tim02} (or \autoref{thm:ex} below).
\end{proof}

The following result is a generalization of Mr\'{o}wka \cite[Theorem 4.11]{Mro68} (see also Blair \cite[Theorem 3.2]{Bla81}), which states that for $f \in C_{b}(S, \mathbb{R})$, $f$ has a continuous extension over $\Omega$ if and only if disjoint Lebesgue sets of $f$ are completely separated in $\Omega$; in particular, if disjoint Lebesgue sets of $f$ are completely separated in $\Omega$, then for every $g \in C_{b}(\Omega, \mathbb{R})$, so are disjoint Lebesgue sets of $f + g$. It would seem that the original proofs due to Mr\'{o}wka and Blair cannot be adapted to the vector-valued case.

\begin{thm}\label{thm:ex}
	Let $f \in C(S, X)$ where $S \subset \Omega$. Suppose that the closure of the linear span of $f(S)$ is metrizable and complete, and $\overline{f(S)}$ is compact. Then $A = C(\Omega, [0, 1])$ separates disjoint Lebesgue sets of $f$ in $S$ (see \autoref{def:sep}) if and only if there is $g \in C(\Omega, X)$ such that $g|_S = f$, $g(\Omega)$ belongs to the convex hull of $f(S)$ and $\overline{g(\Omega)}$ is compact.
\end{thm}

\begin{proof}
	Necessity. This is clear.
	
	Sufficiency. Let $K := \overline{\mathrm{span} f(S)}$ which is metrizable, complete and locally convex according to the hypothesis.
	Let $W_{n}$ be a countable local base of $K$ such that $W_{n + 1} + W_{n + 1} \subset W_{n}$ and $W_{n}$ ($n = 1, 2, ...$) are convex and balanced. In the following, the convex hull of $X_0 \subset K$ is denoted by $co{X_0}$. 
	By \autoref{lem:pre}, there is $\phi_1 \in C_{b}(\Omega, K)$ with $\overline{\phi_1(\Omega)}$ compact such that $\overline{(f - \phi_1) (S)} \subset W_1$ and $\phi_1(\Omega) \subset co f(S)$. For $f_1 = f - \phi_1$, since $\overline{f_1(\Omega)}$ is compact and $A = C(\Omega, [0, 1])$ separates disjoint Lebesgue sets of $f_1$ in $S$, by \autoref{lem:pre}, there is $\phi_2 \in C_{b}(\Omega, K)$ with $\overline{\phi_2(\Omega)}$ compact such that $\overline{(f_1 - \phi_2) (S)} \subset W_2$ and $\phi_2(\Omega) \subset co f_1(S)$. Inductively, for $f_n = f_{n-1} - \phi_{n}$, there is $\phi_{n+1} \in C_{b}(\Omega, K)$ with $\overline{\phi_{n+1}(\Omega)}$ compact such that $\overline{(f_n - \phi_{n+1}) (S)} \subset W_{n+1}$ and $\phi_{n+1}(\Omega) \subset co f_n(S)$, $n = 2, 3, ...$.
	
	Write $u_{n} = \sum_{i = 1}^{n} \phi_i$. Then $\overline{(f - u_{n})(S)} \subset W_{n}$ and $\overline{u_{n}(\Omega)} \subset co ({(f - u_{n-1})(S)}) \subset W_{n-1}$. It yields that $\phi = \lim_{n \to \infty}u_{n}$ which is uniformly convergent in $\Omega$ and $\phi|_S = f$. Since each $\overline{\phi_{n}(\Omega)}$ is compact, $\overline{\phi(\Omega)}$ is compact. Let $K_0 = \overline{f(S)}$. As $K$ is a metric space, by the Dugundji extension theorem (\cite[Chapter IX Theorem 6.1]{Dug66}), there is a continuous function $h: K \to co K_0$ such that $h|_{K_0} = \id$. The function $g = h(\phi)$ satisfies the conclusion. The proof is complete.
\end{proof}

The theory of extension of continuous functions deals with that for which $\Omega$, $S \subset \Omega$, $X$ and $f \in C(S, X)$, there is $g \in C(\Omega, X)$ such that $g|_S = f$; it has been attracted many researchers. We briefly summarize some developments below; see \cite{Hus10} for a more comprehensive overview.
\begin{enumerate}[(a)]
	\item Tietze extension theorem \cite[Chapter VII Theorem 5.1]{Dug66}: $\Omega$ normal, $S$ closed, all $f \in C(S, \mathbb{R})$.
	\item Dugundji extension theorem \cite[Chapter IX Theorem 6.1]{Dug66}: $\Omega$ a metric space, $S$ closed, $X$ a locally convex topological vector space, all $f \in C(S, X)$.
	\item Hanner \cite{Han51}: $\Omega$ normal, $S$ closed, $X$ a \emph{separable} Fr\'{e}chet space, all $f \in C(S, X)$.
	\item Dowker \cite{Dow52}: $\Omega$ collectionwise normal, $S$ closed, $X$ a Fr\'{e}chet space, all $f \in C(S, X)$.
	\item Gillman--Jerison \cite{GJ60}: characterization of $C^*$-embedded subsets, i.e., for which $S$ so that all functions in $ C_{b}(S, \mathbb{R})$ can be extended over $\Omega$ (see also \cite{BH74} for more developments); \autoref{thm:gj} shows that for each $f \in C_{b}(S, X)$ with relatively compact range, $S$ $C^*$-embedded, and $X$ a Fr\'{e}chet space, such extension holds.
	\item Gillman--Jerison \cite{GJ60}: characterization of $C$-embedded subsets, i.e., for which $S$ so that all functions in $ C(S, \mathbb{R})$ can be extended over $\Omega$ (see also \cite{BH74} for more developments); see \cite[Theorem 4.6]{GOY06} for the results concerning the extension of vector-valued functions in this case.
	\item Mr\'{o}wka \cite{Mro68} and Blair \cite{Bla81}: characterization of single $f \in C_{b}(S, \mathbb{R})$ (and $C(S, \mathbb{R})$) which can be extended; \autoref{thm:ex} generalizes it to vector-valued functions having relatively compact ranges.
	\item Hern\'{a}ndez-Mu\~{n}oz \cite[Theorem 9]{Her94}: characterization of single $f \in C(S, X)$ which can be extended where $X$ is a separable Banach space (or Fr\'{e}chet space).
\end{enumerate}

The necessary and sufficient condition on $f \in C(S, X)$ with relatively compact range so that $f$ is extendable is usually weaker than on $f \in C(S, X)$ (without relatively compact range).

\begin{appendices}
	\setcounter{equation}{0}
	\renewcommand{\theequation}{\Alph{section}.\arabic{equation}}
	
\end{appendices}

%\bibliographystyle{amsalpha}



\begin{bibdiv}
\begin{biblist}

\bib{ACIT18}{incollection}{
      author={Angoa, J., Amador},
      author={Contreras, A., Carreto},
      author={Ibarra, M., Contreras},
      author={Tamariz, A., Mascar\'{u}a},
       title={Basic and classic results on pseudocompact spaces},
        date={2018},
   booktitle={Pseudocompact topological spaces},
      series={Dev. Math.},
      volume={55},
   publisher={Springer, Cham},
       pages={1\ndash 38},
      review={\MR{3822416}},
}

\bib{Are49}{article}{
      author={Arens, Richard},
       title={Approximation in, and representation of, certain {B}anach
  algebras},
        date={1949},
        ISSN={0002-9327},
     journal={Amer. J. Math.},
      volume={71},
       pages={763\ndash 790},
         url={https://doi.org/10.2307/2372363},
      review={\MR{32953}},
}

\bib{BD81}{article}{
      author={Brosowski, Bruno},
      author={Deutsch, Frank},
       title={An elementary proof of the {S}tone-{W}eierstrass theorem},
        date={1981},
        ISSN={0002-9939},
     journal={Proc. Amer. Math. Soc.},
      volume={81},
      number={1},
       pages={89\ndash 92},
         url={https://doi.org/10.2307/2043993},
      review={\MR{589143}},
}

\bib{BH74}{article}{
      author={Blair, Robert~L.},
      author={Hager, Anthony~W.},
       title={Extensions of zero-sets and of real-valued functions},
        date={1974},
        ISSN={0025-5874,1432-1823},
     journal={Math. Z.},
      volume={136},
       pages={41\ndash 52},
         url={https://doi.org/10.1007/BF01189255},
      review={\MR{385793}},
}

\bib{Bie73}{article}{
      author={Bierstedt, Klaus-Dieter},
       title={Gewichtete {R}\"{a}ume stetiger vektorwertiger {F}unktionen und
  das injektive {T}ensorprodukt. {I}},
        date={1973},
        ISSN={0075-4102,1435-5345},
     journal={J. Reine Angew. Math.},
      volume={259},
       pages={186\ndash 210},
         url={https://doi.org/10.1515/crll.1973.259.186},
      review={\MR{318871}},
}

\bib{Bie81}{incollection}{
      author={Bierstedt, Klaus-D.},
       title={The approximation-theoretic localization of {S}chwartz's
  approximation property for weighted locally convex function spaces and some
  examples},
        date={1981},
   booktitle={Functional analysis, holomorphy, and approximation theory
  ({P}roc. {S}em., {U}niv. {F}ed. {R}io de {J}aneiro, {R}io de {J}aneiro,
  1978)},
      series={Lecture Notes in Math.},
      volume={843},
   publisher={Springer, Berlin},
       pages={93\ndash 149},
      review={\MR{610827}},
}

\bib{Bis61}{article}{
      author={Bishop, Errett},
       title={A generalization of the {S}tone-{W}eierstrass theorem},
        date={1961},
        ISSN={0030-8730},
     journal={Pacific J. Math.},
      volume={11},
       pages={777\ndash 783},
         url={http://projecteuclid.org/euclid.pjm/1103037116},
      review={\MR{0133676}},
}

\bib{Bla81}{article}{
      author={Blair, Robert~L.},
       title={Extensions of {L}ebesgue sets and of real-valued functions},
        date={1981},
        ISSN={0011-4642},
     journal={Czechoslovak Math. J.},
      volume={31(106)},
      number={1},
       pages={63\ndash 74},
        note={With a loose Russian summary},
      review={\MR{604112}},
}

\bib{BM00}{article}{
      author={Bustamante, Jorge},
      author={Montalvo, Francisco},
       title={Stone-{W}eierstrass theorems in {$C^\ast(X)$}},
        date={2000},
        ISSN={0021-9045},
     journal={J. Approx. Theory},
      volume={107},
      number={1},
       pages={143\ndash 159},
         url={https://doi.org/10.1006/jath.2000.3507},
      review={\MR{1799556}},
}

\bib{BM83}{article}{
      author={Blasco, J.~L.},
      author={Molt\'{o}, A.},
       title={On the uniform closure of a linear space of bounded real-valued
  functions},
        date={1983},
        ISSN={0003-4622},
     journal={Ann. Mat. Pura Appl. (4)},
      volume={134},
       pages={233\ndash 239},
         url={https://doi.org/10.1007/BF01773506},
      review={\MR{736741}},
}

\bib{BP20}{book}{
      author={Bucur, Ileana},
      author={Paltineanu, Gavriil},
       title={Topics in uniform approximation of continuous functions},
      series={Frontiers in Mathematics},
   publisher={Birkh\"{a}user/Springer, Cham},
        date={2020},
        ISBN={978-3-030-48411-8; 978-3-030-48412-5},
         url={https://doi.org/10.1007/978-3-030-48412-5},
      review={\MR{4311187}},
}

\bib{Buc58}{article}{
      author={Buck, R.~Creighton},
       title={Bounded continuous functions on a locally compact space},
        date={1958},
        ISSN={0026-2285},
     journal={Michigan Math. J.},
      volume={5},
       pages={95\ndash 104},
         url={http://projecteuclid.org/euclid.mmj/1028998054},
      review={\MR{0105611}},
}

\bib{Che18f}{article}{
      author={Chen, Deliang},
       title={A note on {M}achado-{B}ishop theorem in weighted spaces with
  applications},
        date={2019},
        ISSN={0021-9045},
     journal={J. Approx. Theory},
      volume={247},
       pages={1\ndash 19},
         url={https://doi.org/10.1016/j.jat.2019.07.004},
      review={\MR{3984166}},
}

\bib{dBra59}{article}{
      author={de~Branges, Louis},
       title={The {S}tone-{W}eierstrass theorem},
        date={1959},
        ISSN={0002-9939,1088-6826},
     journal={Proc. Amer. Math. Soc.},
      volume={10},
       pages={822\ndash 824},
         url={https://doi.org/10.2307/2033481},
      review={\MR{113131}},
}

\bib{Dow52}{article}{
      author={Dowker, C.~H.},
       title={On a theorem of {H}anner},
        date={1952},
        ISSN={0004-2080},
     journal={Ark. Mat.},
      volume={2},
       pages={307\ndash 313},
         url={https://doi.org/10.1007/BF02591500},
      review={\MR{50874}},
}

\bib{Dug66}{book}{
      author={Dugundji, James},
       title={Topology},
   publisher={Allyn and Bacon, Inc., Boston, Mass.},
        date={1966},
      review={\MR{0193606}},
}

\bib{GJ60}{book}{
      author={Gillman, Leonard},
      author={Jerison, Meyer},
       title={Rings of continuous functions},
      series={The University Series in Higher Mathematics},
   publisher={D. Van Nostrand Co., Inc., Princeton, N.J.-Toronto-London-New
  York},
        date={1960},
      review={\MR{0116199}},
}

\bib{Gli59}{article}{
      author={Glicksberg, Irving},
       title={Stone-\v{C}ech compactifications of products},
        date={1959},
        ISSN={0002-9947,1088-6850},
     journal={Trans. Amer. Math. Soc.},
      volume={90},
       pages={369\ndash 382},
         url={https://doi.org/10.2307/1993177},
      review={\MR{105667}},
}

\bib{Gli63}{article}{
      author={Glicksberg, I.},
       title={Bishop's generalized {S}tone-{W}eierstrass theorem for the strict
  topology},
        date={1963},
        ISSN={0002-9939},
     journal={Proc. Amer. Math. Soc.},
      volume={14},
       pages={329\ndash 333},
         url={https://doi.org/10.2307/2034636},
      review={\MR{0146645}},
}

\bib{GM93}{article}{
      author={Garrido, M.~I.},
      author={Montalvo, F.},
       title={On some generalizations of the {K}akutani-{S}tone and
  {S}tone-{W}eierstrass theorems},
        date={1993},
        ISSN={0236-5294},
     journal={Acta Math. Hungar.},
      volume={62},
      number={3-4},
       pages={199\ndash 208},
         url={https://doi.org/10.1007/BF01874642},
      review={\MR{1250902}},
}

\bib{GOY06}{article}{
      author={Gutev, Valentin},
      author={Ohta, Haruto},
      author={Yamazaki, Kaori},
       title={Extensions by means of expansions and selections},
        date={2006},
        ISSN={0927-6947},
     journal={Set-Valued Anal.},
      volume={14},
      number={1},
       pages={69\ndash 104},
         url={https://doi.org/10.1007/s11228-005-0008-y},
      review={\MR{2232459}},
}

\bib{Hag66}{article}{
      author={Hager, Anthony~W.},
       title={Some remarks on the tensor product of function rings},
        date={1966},
        ISSN={0025-5874},
     journal={Math. Z.},
      volume={92},
       pages={210\ndash 224},
         url={https://doi.org/10.1007/BF01111186},
      review={\MR{193613}},
}

\bib{Han51}{article}{
      author={Hanner, Olof},
       title={Solid spaces and absolute retracts},
        date={1951},
        ISSN={0004-2080},
     journal={Ark. Mat.},
      volume={1},
       pages={375\ndash 382},
         url={https://doi.org/10.1007/BF02591374},
      review={\MR{43458}},
}

\bib{Hew47}{article}{
      author={Hewitt, Edwin},
       title={Certain generalizations of the {W}eierstrass approximation
  theorem},
        date={1947},
        ISSN={0012-7094},
     journal={Duke Math. J.},
      volume={14},
       pages={419\ndash 427},
         url={http://projecteuclid.org/euclid.dmj/1077474139},
      review={\MR{21662}},
}

\bib{Her94}{article}{
      author={Hern\'{a}ndez-Mu\~{n}oz, Salvador},
       title={Approximation and extension of continuous functions},
        date={1994},
        ISSN={0263-6115},
     journal={J. Austral. Math. Soc. Ser. A},
      volume={57},
      number={2},
       pages={149\ndash 157},
      review={\MR{1288670}},
}

\bib{Hus10}{article}{
      author={Hu\v{s}ek, M.},
       title={Extension of mappings and pseudometrics},
        date={2010},
        ISSN={0213-8743},
     journal={Extracta Math.},
      volume={25},
      number={3},
       pages={277\ndash 308},
      review={\MR{2857999}},
}

\bib{KG06}{incollection}{
      author={Kundu, S.},
      author={Garg, Pratibha},
       title={The pseudocompact-open topology on {$C(X)$}},
        date={2006},
      volume={30},
       pages={279\ndash 299},
        note={Spring Topology and Dynamical Systems Conference},
      review={\MR{2280673}},
}

\bib{Kle75}{article}{
      author={Kleinst\"{u}ck, Gert},
       title={Der beschr\"{a}nkte {F}all des gewichteten {A}pproximations
  problems f\"{u}r vektorwertige {F}unktionen},
        date={1975},
        ISSN={0025-2611,1432-1785},
     journal={Manuscripta Math.},
      volume={17},
      number={2},
       pages={123\ndash 149},
         url={https://doi.org/10.1007/BF01154086},
      review={\MR{407592}},
}

\bib{Mac77}{article}{
      author={Machado, Silvio},
       title={On {B}ishop's generalization of the {W}eierstrass-{S}tone
  theorem},
        date={1977},
     journal={Indag. Math.},
      volume={39},
      number={3},
       pages={218\ndash 224},
      review={\MR{0448046}},
}

\bib{Mey67}{article}{
      author={Meyer, Paul~R.},
       title={Topologies with the {S}tone-{W}eierstrass property},
        date={1967},
        ISSN={0002-9947},
     journal={Trans. Amer. Math. Soc.},
      volume={126},
       pages={236\ndash 243},
         url={https://doi.org/10.2307/1994450},
      review={\MR{203699}},
}

\bib{Mro68}{article}{
      author={Mr\'{o}wka, S.},
       title={On some approximation theorems},
        date={1968},
        ISSN={0028-9825},
     journal={Nieuw Arch. Wisk. (3)},
      volume={16},
       pages={94\ndash 111},
      review={\MR{244938}},
}

\bib{Nac65}{article}{
      author={Nachbin, Leopoldo},
       title={Weighted approximation for algebras and modules of continuous
  functions: {R}eal and self-adjoint complex cases},
        date={1965},
        ISSN={0003-486X},
     journal={Ann. of Math. (2)},
      volume={81},
       pages={289\ndash 302},
         url={https://doi.org/10.2307/1970617},
      review={\MR{0176353}},
}

\bib{Nel68}{article}{
      author={Nel, L.~D.},
       title={Theorems of {S}tone-{W}eierstrass type for non-compact spaces},
        date={1968},
        ISSN={0025-5874},
     journal={Math. Z.},
      volume={104},
       pages={226\ndash 230},
         url={https://doi.org/10.1007/BF01110290},
      review={\MR{226384}},
}

\bib{Pczy57}{article}{
      author={Pe\l~czy\'{n}ski, A.},
       title={A generalisation of {S}tone's theorem on approximation},
        date={1957},
     journal={Bull. Acad. Polon. Sci. Cl. III.},
      volume={5},
       pages={105\ndash 107, X},
      review={\MR{86166}},
}

\bib{PM73}{article}{
      author={Prolla, Jo\~{a}o~Bosco},
      author={Machado, Silvio},
       title={Weighted {G}rothendieck subspaces},
        date={1973},
        ISSN={0002-9947},
     journal={Trans. Amer. Math. Soc.},
      volume={186},
       pages={247\ndash 258 (1974)},
         url={https://doi.org/10.2307/1996564},
      review={\MR{0402477}},
}

\bib{Pro71}{article}{
      author={Prolla, Jo\~{a}o~B.},
       title={Bishop's generalized {S}tone-{W}eierstrass theorem for weighted
  spaces},
        date={1971},
        ISSN={0025-5831},
     journal={Math. Ann.},
      volume={191},
       pages={283\ndash 289},
         url={https://doi.org/10.1007/BF01350331},
      review={\MR{0290015}},
}

\bib{Pro77}{book}{
      author={Prolla, Jo\~{a}o~Bosco},
       title={Approximation of vector valued functions},
   publisher={North-Holland Publishing Co., Amsterdam-New York-Oxford},
        date={1977},
        ISBN={0-444-85030-9},
        note={North-Holland Mathematics Studies, Vol. 25, Notas de
  Matem\'{a}tica, No. 61. [Notes on Mathematics, No. 61]},
      review={\MR{0500122}},
}

\bib{Pro88}{article}{
      author={Prolla, Jo\~{a}o~B.},
       title={A generalized {B}ernstein approximation theorem},
        date={1988},
        ISSN={0305-0041},
     journal={Math. Proc. Cambridge Philos. Soc.},
      volume={104},
      number={2},
       pages={317\ndash 330},
         url={https://doi.org/10.1017/S030500410006549X},
      review={\MR{948917}},
}

\bib{Pro94}{article}{
      author={Prolla, Jo\~{a}o~B.},
       title={On the {W}eierstrass-{S}tone theorem},
        date={1994},
        ISSN={0021-9045},
     journal={J. Approx. Theory},
      volume={78},
      number={3},
       pages={299\ndash 313},
         url={https://doi.org/10.1006/jath.1994.1080},
      review={\MR{1292963}},
}

\bib{Ran84}{article}{
      author={Ransford, T.~J.},
       title={A short elementary proof of the {B}ishop-{S}tone-{W}eierstrass
  theorem},
        date={1984},
        ISSN={0305-0041},
     journal={Math. Proc. Cambridge Philos. Soc.},
      volume={96},
      number={2},
       pages={309\ndash 311},
         url={https://doi.org/10.1017/S0305004100062204},
      review={\MR{757664}},
}

\bib{Rud91}{book}{
      author={Rudin, Walter},
       title={Functional analysis},
     edition={Second},
      series={International Series in Pure and Applied Mathematics},
   publisher={McGraw-Hill, Inc., New York},
        date={1991},
        ISBN={0-07-054236-8},
      review={\MR{1157815}},
}

\bib{San18}{incollection}{
      author={Sanchis, M.},
       title={Bounded subsets of {T}ychonoff spaces: a survey of results and
  problems},
        date={2018},
   booktitle={Pseudocompact topological spaces},
      series={Dev. Math.},
      volume={55},
   publisher={Springer, Cham},
       pages={107\ndash 150},
      review={\MR{3822419}},
}

\bib{Sum71}{article}{
      author={Summers, W.~H.},
       title={The general complex bounded case of the strict weighted
  approximation problem},
        date={1971},
        ISSN={0025-5831},
     journal={Math. Ann.},
      volume={192},
       pages={90\ndash 98},
         url={https://doi.org/10.1007/BF02052753},
      review={\MR{0284800}},
}

\bib{Tam61}{article}{
      author={Tamano, Hisahiro},
       title={A note on the pseudo-compactness of the product of two spaces},
        date={1960/61},
        ISSN={0368-8887},
     journal={Mem. Coll. Sci. Univ. Kyoto Ser. A. Math.},
      volume={33},
       pages={225\ndash 230},
         url={https://doi.org/10.1215/kjm/1250775908},
      review={\MR{120619}},
}

\bib{Tim02}{article}{
      author={Timofte, Vlad},
       title={Special uniform approximations of continuous vector-valued
  functions. {I}. {S}pecial approximations in {$C_X(T)$}},
        date={2002},
        ISSN={0021-9045,1096-0430},
     journal={J. Approx. Theory},
      volume={119},
      number={2},
       pages={291\ndash 299},
         url={https://doi.org/10.1006/jath.2002.3733},
      review={\MR{1939286}},
}

\bib{Tim03}{article}{
      author={Timofte, Vlad},
       title={Special uniform approximations of continuous vector-valued
  functions. {II}. {S}pecial approximations in {$C_X(T)\otimes C_Y(S)$}},
        date={2003},
        ISSN={0021-9045,1096-0430},
     journal={J. Approx. Theory},
      volume={123},
      number={2},
       pages={270\ndash 275},
         url={https://doi.org/10.1016/S0021-9045(03)00103-5},
      review={\MR{1990100}},
}

\bib{Tre67}{book}{
      author={Tr\`eves, Fran\c{c}ois},
       title={Topological vector spaces, distributions and kernels},
   publisher={Academic Press, New York-London},
        date={1967},
      review={\MR{225131}},
}

\end{biblist}
\end{bibdiv}
\end{document}